\documentclass[a4paper,11pt]{article}
\usepackage[utf8]{inputenc}

\usepackage{bm}
\usepackage{amsmath}
\usepackage{amsthm}

\usepackage{mathrsfs}
\usepackage{graphicx}
\usepackage{epsfig, setspace}
\usepackage{pst-plot}
\usepackage{booktabs}
\usepackage{multirow}
\usepackage{lscape}
\usepackage{subfig}
\usepackage{url}
\usepackage{setspace}
\usepackage{graphicx, transparent, color}
\usepackage[ruled]{algorithm2e}
\usepackage{etex}

\DeclareSymbolFont{rsfs}{U}{rsfs}{m}{n}
\DeclareSymbolFontAlphabet{\mathcal}{rsfs}

\theoremstyle{definition}
\theoremstyle{plain}
\usepackage[margin=2.5cm]{geometry}
\usepackage{tikz}
\usetikzlibrary{shapes,calc}
\usepackage{verbatim}
\usepackage{lineno}
\usepackage{listings}
\usepackage{textcomp}

\usepackage{epsfig}
\usepackage{amsmath}
\usepackage{pst-plot}
\usepackage{amssymb}
\usepackage{amsfonts}
\usepackage{lmodern}
\usepackage{euscript}
\usepackage{eufrak}
\usepackage{oldgerm}

\numberwithin{equation}{section}

\newtheorem{example}{Example}[section]
\newtheorem{definition}{Definition}[section]

\usepackage[nottoc]{tocbibind}
\usepackage{color}
\pagecolor{white}
\usepackage{slantsc}
\begin{document}

\title{Adaptive IQ and IMQ-RBFs for solving Initial Value Problems: Adam-Bashforth and Adam-Moulton methods}
 \date{}
\author{
{ Samala Rathan\thanks{Department of Humanities and Sciences, Indian Institute of Petroleum and Energy, Visakhapatnam, India-530003 (Email: rathans.math@iipe.ac.in)}, 
Deepit Shah\thanks{Department of Petroleum Engineering and Earth Sciences, Indian Institute of Petroleum and Energy, Visakhapatnam, India-530003}, 
T. Hemanth Kumar\thanks{Department of Chemical Engineering, Indian Institute of Petroleum and Energy, Visakhapatnam, India-530003},
K. Sandeep Charan\thanks{Department of Chemical Engineering, National Institute of Technology, Tiruchirappalli, India-620015}
}}
\maketitle
\begin{abstract}
In this paper, our objective is primarily to use adaptive inverse-quadratic (IQ) and inverse-multi-quadratic (IMQ) radial basis function (RBF) interpolation techniques to develop an enhanced Adam-Bashforth and Adam-Moulton methods. By utilizing a free parameter involved in the radial basis function, the local convergence of the numerical solution is enhanced by making the local truncation error vanish. Consistency and stability analysis is presented along with some numerical results to back up our assertions. The accuracy and rate of convergence of each proposed technique are equal to or better than the original Adam-Bashforth and Adam-Moulton methods by eliminating the local truncation error thus, the proposed adaptive methods are optimal. We conclude that both IQ and IMQ-RBF methods yield an improved order of convergence than classical methods, while the superiority of one method depends on the method and the problem considered.
\end{abstract}
{\small \textbf{Keywords:} Finite difference method, Radial basis interpolation, Multistep method, Stability, Order of accuracy, Rate of convergence..\\[0.5pt]\\
\textbf{AMS subject classification:} 41A10, 65L05, 65L12, 65L20 }

\section{Introduction}
\label{sec:intro}
In this work, we propose  inverse quadratic (IQ) and inverse multi-quadratic (IMQ) RBF interpolation methods to solve first order IVPs. The classical finite difference methods cannot use the local information of solutions to increase the accuracy of solution. To overcome this, the use of RBF based interpolation schemes were explored in literature. The success of RBF interpolation techniques to enhance the order of accuracy of the numerical solution schemes already demonstrated. For example, in \cite{guo2017rbf, guo2017radial} RBF interpolation was used to improve the accuracy of ENO and WENO schemes for solving hyperbolic partial differential equations. In \cite{rathanshah2022,gu2021adaptive,gu2020adaptive} RBF interpolation  is used to improve the accuracy of the finite difference schemes to solve the first order IVPs.

The key principle behind techniques modified with adaptive RBF solvers is the radial basis functions (RBF). In RBF interpolation, the free shape parameter is exploited and allowed to change its value based on the local conditions of the solution. These changes in free shape parameter help to improve the accuracy. The optimal values of the free shape parameter will be obtained by solving an optimization problem. The objective function of an optimization problem is to minimize the leading truncation error. The free shape parameter is the decision variable. The local information on solution is used to define the objective function.

The use of RBF interpolation finite difference schemes to solve the first order IVPs were first proposed in \cite{gu2020adaptive}. Multi quadratic RBF interpolation was used to modify the finite difference methods to improve the accuracy of solutions. The methods like Euler, midpoint, Adams-Bashforth (AB) and Adams-Moultan (AM) methods were developed using the RBF interpolation. The developed methods using the mutli-quadratic RBF interpolation were reported to have higher accuracy compared to their polynomial expansion counter parts. To supplement their previous studies, Gu and Jung \cite{gu2021adaptive} have proposed the use of Gaussian RBF interpolation to derive Euler, midpoint, Adams-Bashforth and Adams-Moultan methods for solving first order IVP problems. For both multi-quadratic RBF and Gaussian RBF methods, an increased accuracy was reported compared to the regular polynomial function based methods. In \cite{rathanshah2022}, second order time stepping methods  with the adaptive inverse quadratic and inverse multi quadratic radial basis function interpolation technique for solving IVPs were studied. The consistency, stability and convergence analysis was also elaborated with its advantages depending on the considered problem and method. Here, we study an extension to the literature available on the application of RBFs to solve the IVPs given in \cite{rathanshah2022} for higher order numerical methods. 

 In this work, we have used inverse quadratic and inverse multi-quadratic RBF interpolation methods to solve first order IVPs. The classical finite difference techniques (AB-2, AB-3, AM-2 and AM-3) to solve IVPs were developed with the proposed RBFs and the effect of this development on the accuracy of the solution was studied.  The shape parameter was locally optimised to make the leading error term(s) disappear. This will lower the local truncation error and subsequently reduce the global error. The optimal free shape parameter obtained here increased the order of convergence under the assumption of smooth solution. Further, the stability areas of proposed and original techniques were compared. We limit our consideration to the scenarios where the RBFs in consideration have just one shape parameter, while several shape parameters may be used for further improvement. The polynomial interpolation is a limit instance of the RBF interpolation, which is one of the features of the proposed techniques using RBFs. That is, if the shape parameter disappears or does not fit, based on the RBF definition used \cite{rathanshah2022,gu2021adaptive,gu2020adaptive}, the RBF interpolation becomes equal to the polynomial interpolation. This indicates that even if the shape parameter does not match or deviates from the ideal value, the proposed techniques demonstrate at least the same rate of convergence as the original methods. As a result, with the proposed schemes, the convergence characteristics can only get better. 

The following is a breakdown of the paper's structure. The RBF interpolation is briefly discussed in Section 2. In Section 3, we describe proposed RBF interpolation-based techniques for solving IVPs. We also describe how the resulting techniques are consistent and stable in Section 4.  An overview of the proposed techniques, and comparative study with the actual methods can also be seen in Section 4. We demonstrate several numerical experiments in Section 5.
Section 6 discusses a brief conclusion and future work. 

\section{RBF Interpolation}
Let us consider the construction of radial basis function(RBF) interpolation in one-space dimension. Given $(N+1)$ distinct data points $(x_0, u_0), . . . ,  (x_N, u_N )$  with $u_k$ the value of the unknown function $u(x)$ at $x = x_k$, where $ x \in \mathbb{R}$. We use the RBFs, $\phi_k (x) = \phi \big(|x - x_k|, \epsilon_k \big)$, where  $\epsilon_k$ is a shape parameter, to find an interpolant based on the given $(N+1)$ data points. The value of $\epsilon_k$ can vary over $x_k$. The interpolant $r (x)$ takes the form of a weighted sum of RBFs
\begin{equation}\label{eq1}
r(x)=\sum_{k=0}^{N} \lambda_k\phi \big(|x - x_k |, \epsilon_k \big)
\end{equation}
where $\lambda_k$ are the unknown expansion parameters to be determined.
Using interpolation restraints $$r(x_k)=u_k, k= 0,1,...,N,$$ the expansion coefficients $\lambda_k$ satisfy the following linear system
\begin{eqnarray*}
\begin{pmatrix}
\phi(|x_0 - x_0|, \epsilon_0) & \phi( |x_0 - x_1|, \epsilon_1) &.... & \phi(|x_0 - x_N |, \epsilon_N )\\
\phi(|x_1 - x_0|, \epsilon_0) & \phi( |x_1 - x_1|, \epsilon_1) &....& \phi(|x_1 - x_N |, \epsilon_N )\\
... & ... &.... &....\\
... & ... &.... &....\\
... & ... &.... &....\\
\phi(|x_N - x_0|, \epsilon_0) &  \phi( |x_N - x_1|, \epsilon_1) &....& \phi(|x_N - x_N |, \epsilon_N )
\end{pmatrix}
\begin{pmatrix}
\lambda_0\\\lambda_1\\.\\.\\.\\\lambda_N
\end{pmatrix}=
\begin{pmatrix}
u_0\\u_1\\.\\.\\.\\u_N\end{pmatrix}
\end{eqnarray*}
Until and unless specified explicitly, we consider all the shape parameters $\epsilon_k$'s are same, i.e., $\epsilon_k=\epsilon$, for all $k$. In this paper, we use inverse multi-quadratic (IMQ) and inverse quadratic (IQ)-RBFs for solving initial value problems.

\subsection{IMQ-RBF Interpolation}
Let us consider the inverse-multi-quadratic (IMQ) radial basis function $\phi_k(x)=\displaystyle\frac{1}{\sqrt{1+\epsilon_k^2(x-x_k)^2}}$. Now, we derive the interpolation for the $N=1$ and $N=2$ data points.
\subsubsection{N=1}
Consider the interpolant for $N=1$, we have
\begin{equation}\label{eq2}
r(x)=\displaystyle\lambda_0\phi_0(x)+\lambda_1\phi_1(x).
\end{equation}
Using the interpolation condition $r(x_k)=u_k, k= 0,1,$ the interpolation matrix becomes a symmetric matrix with all diagonal entries  $1$.  Thus, we have
\begin{equation*}
\displaystyle\begin{pmatrix}
1&\displaystyle\frac{1}{\sqrt{1+\epsilon^2h^2}}\\\displaystyle\frac{1}{\sqrt{1+\epsilon^2h^2}}&1
\end{pmatrix}
\begin{pmatrix}
\displaystyle\lambda_0\\\lambda_1
\end{pmatrix}=\displaystyle\begin{pmatrix}
u_0\\u_1
\end{pmatrix},
\end{equation*}
where $h=x_1-x_0.$  Solving for $\lambda_k, k=0,1$, we get
$$\displaystyle\lambda_0=\frac{1+\epsilon^2h^2}{\epsilon^2h^2}\bigg(u_0-\frac{u_1}{\sqrt{1+\epsilon^2h^2}}\bigg),$$
$$\displaystyle\lambda_1=\frac{1+\epsilon^2h^2}{\epsilon^2h^2}\bigg(u_1-\frac{u_0}{\sqrt{1+\epsilon^2h^2}}\bigg).$$
Differentiating the interpolant $r(x)$  with respect to $x$, we get
\begin{equation*}\label{eq3}
\frac{d}{dx}r(x)=-\frac{\lambda_0(x-x_0)\epsilon^2}{{(1+\epsilon^2(x-x_0)^2})^\frac{3}{2}}-\frac{\lambda_1(x-x_1)\epsilon^2}{{(1+\epsilon^2(x-x_1)^2})^\frac{3}{2}},
\end{equation*}
and evaluating at $x=x_0$, we have
\begin{equation}\label{eq4}
\frac{d}{dx} r(x_0)=\frac{\lambda_1\epsilon^2h}{{(1+\epsilon^2h^2})^\frac{3}{2}}.
\end{equation}
Substituting the value of $\lambda_1$ in \eqref{eq4}, we have
\begin{equation}\label{eqimq}
\frac{d}{dx} r(x_0)=\frac{u_1\sqrt{1+\epsilon^2h^2}-u_0}{(1+\epsilon^2h^2)h},
\end{equation}
and as $\epsilon \to 0$, the equation \eqref{eqimq} become
\begin{equation*}
\lim_{\epsilon \to 0}\frac{d}{dx}r(x_0)=\frac{u_1-u_0}{h}.
\end{equation*}
Thus, the method reduces to forward difference formula of $u^{\prime}_x$ at $x=x_0.$\\

\subsubsection{N=2}
Consider the interpolant for $N=2$, we get
\begin{equation}\label{eq5}
r(x)=\displaystyle\lambda_0\phi_0(x)+\lambda_1\phi_1(x)+\lambda_2\phi_2(x).
\end{equation}
Using the interpolation condition $r(x_k)=u_k, k= 0,1,2$ the interpolation matrix becomes a symmetric matrix with all diagonal entries  $1$.  Thus, we have
\begin{equation*}
\displaystyle\begin{pmatrix}
1&\displaystyle\frac{1}{\sqrt{1+\epsilon^2h^2}}&\displaystyle\frac{1}{\sqrt{1+4\epsilon^2h^2}}\\\displaystyle\frac{1}{\sqrt{1+\epsilon^2h^2}}&1&\displaystyle\frac{1}{\sqrt{1+\epsilon^2h^2}}\\\displaystyle\frac{1}{\sqrt{1+4\epsilon^2h^2}}&\displaystyle\frac{1}{\sqrt{1+\epsilon^2h^2}}&1
\end{pmatrix}
\begin{pmatrix}
\displaystyle\lambda_0\\\lambda_1\\\lambda_2
\end{pmatrix}=\displaystyle\begin{pmatrix}
u_0\\u_1\\u_2
\end{pmatrix},
\end{equation*}
where $h=x_{k+1}-x_k.$  Solving for $\lambda_k, k=0,1,2$, we get
\begin{eqnarray*}
\begin{split}
 \lambda_0 &=\frac{1+4\epsilon^2h^2}{\big(2\sqrt{1+\epsilon^2h^2}\big)\big(1-\sqrt{1+4\epsilon^2h^2}+2\epsilon^4h^4\sqrt{1+4\epsilon^2h^2}-2\epsilon^2h^2
  \big(-2+\sqrt{1+4\epsilon^2h^2}\big)\big)}\bigg(\epsilon^2h^2\sqrt{1+\epsilon^2h^2}\\
  & \sqrt{1+4\epsilon^2h^2}u_0-\big(1+\epsilon^2h^2\big)\big(-1+\sqrt{1+4\epsilon^2h^2}\big)u_1+\sqrt{1+\epsilon^2h^2}\big(-1-\epsilon^2h^2+\sqrt{1+4\epsilon^2h^2}\big)u_2\bigg),\\
\lambda_1 &= \frac{1}{\big(2-2\sqrt{1+4\epsilon^2h^2}+4\epsilon^4h^4\sqrt{1+4\epsilon^2h^2}
-4\epsilon^2h^2\big(-2+\sqrt{1+4\epsilon^2h^2}\big)\big)}\bigg(\sqrt{1+\epsilon^2h^2}\big(-\big(\big(1+4\epsilon^2h^2\big)\\
& \big(-1+\sqrt{1+4\epsilon^2h^2}\big)u_0\big)
+4\epsilon^2h^2\sqrt{1+\epsilon^2h^2}\sqrt{1+4\epsilon^2h^2}u_1-\big(1+4\epsilon^2h^2\big)\big(-1+\sqrt{1+4\epsilon^2h^2}\big)u_2\big)\bigg),\\
 \lambda_2 &= \frac{1+4\epsilon^2h^2}{\big(2\sqrt{1+\epsilon^2h^2}\big)\big(1-\sqrt{1+4\epsilon^2h^2}+2\epsilon^4h^4\sqrt{1+4\epsilon^2h^2}-2\epsilon^2h^2
\big(-2+\sqrt{1+4\epsilon^2h^2}\big)\big)}\bigg(\sqrt{1+\epsilon^2h^2}\\
& \big(-1-\epsilon^2h^2+\sqrt{1+4\epsilon^2h^2}\big)u_0-\big(1+\epsilon^2h^2\big)\big(-1+\sqrt{1+4\epsilon^2h^2}\big)u_1+\epsilon^2h^2\sqrt{1+\epsilon^2h^2}\sqrt{1+4\epsilon^2h^2}u_2\bigg).
\end{split}
\end{eqnarray*}
Differentiating the interpolant $r(x)$  with respect to $x$, we get
\begin{equation*}\label{eq6}
\frac{d}{dx}r(x)=-\frac{\lambda_0(x-x_0)\epsilon^2}{{(1+\epsilon^2(x-x_0)^2})^\frac{3}{2}}-\frac{\lambda_1(x-x_1)\epsilon^2}{{(1+\epsilon^2(x-x_1)^2})^\frac{3}{2}}-\frac{\lambda_2(x-x_2)\epsilon^2}{{(1+\epsilon^2(x-x_2)^2})^\frac{3}{2}},
\end{equation*}
and evaluating at $x=x_1$, we have
\begin{equation}\label{eq7}
\frac{d}{dx} r(x_1)=\frac{\epsilon^2h}{{(1+\epsilon^2h^2})^\frac{3}{2}}(\lambda_2-\lambda_0).
\end{equation}
Substituting the value of $\lambda_0,\lambda_2$ in \eqref{eq7}, we have
\begin{equation}\label{eqimq2}
\frac{d}{dx} r(x_1)=\frac{\big(1+4\epsilon^2h^2+\sqrt{1+4\epsilon^2h^2}\big)}{4h\big(1+\epsilon^2h^2\big)^{3/2}}\bigg(u_2-u_0\bigg),
\end{equation}
and as $\epsilon \to 0$, the equation \eqref{eqimq2} become
\begin{equation*}
\lim_{\epsilon \to 0}\frac{d}{dx}r(x_1)=\frac{u_2-u_0}{2h}.
\end{equation*}
Thus, the method reduces to centred difference formula of $u^{\prime}_x$ at $x=x_1.$

\subsection{IQ-RBF Interpolation}
The inverse-quadratic (IQ) radial basis function is defined as $\phi_k(x)=\displaystyle\frac{1}{1+\epsilon_k^2(x-x_k)^2}$. Here, we do the interpolation for the case of $N=1$ and $N=2$.
\subsubsection{N=1}
Consider the interpolant ,
 \begin{equation}\label{eq8}
r(x)=\displaystyle\lambda_0\phi_0(x)+\lambda_1\phi_1(x),
\end{equation}
the interpolation matrix becomes a symmetric matrix with all diagonal entries 1, which is
\begin{equation*}
\begin{pmatrix}
1&\displaystyle\frac{1}{1+\epsilon^2h^2}\\\displaystyle\frac{1}{1+\epsilon^2h^2}&1
\end{pmatrix}
\begin{pmatrix}
\displaystyle\lambda_0\\\lambda_1
\end{pmatrix}=\displaystyle\begin{pmatrix}
u_0\\u_1
\end{pmatrix}.
\end{equation*}
Solving for $\lambda_k, k=0,1$, we get
$$\displaystyle\lambda_0=\frac{(1+\epsilon^2h^2)}{\epsilon^2h^2(2+\epsilon^2h^2)}\bigg((1+\epsilon^2h^2)u_0-u_1\bigg),$$
$$\displaystyle\lambda_1=\frac{(1+\epsilon^2h^2)}{\epsilon^2h^2(2+\epsilon^2h^2)}\bigg((1+\epsilon^2h^2)u_1-u_0\bigg).$$
Differentiating $r(x)$ with respect to $x$, we get
\begin{equation}\label{eq9}
\frac{d}{dx}r(x)=-2\frac{\lambda_0(x-x_0)\epsilon^2}{{(1+\epsilon^2(x-x_0)^2})^2}-2\frac{\lambda_1(x-x_1)\epsilon^2}{{(1+\epsilon^2(x-x_1)^2})^2}.\end{equation}
At $x=x_0$ in \eqref{eq9}, we obtain
\begin{equation}\label{eq10}
\frac{d}{dx} r(x_0)=2\frac{\lambda_1\epsilon^2h}{{(1+\epsilon^2h^2})^2}.\end{equation}
Substituting value of $\lambda_1$ in \eqref{eq10}, we get
\begin{equation}\label{eqiq}
\frac{d}{dx} r(x_0)=\dfrac{2\bigg(u_1(1+\epsilon^2h^2)-u_0\bigg)}{\bigg(1+\epsilon^2h^2\bigg)\bigg(2+\epsilon^2h^2\bigg)h},
\end{equation}
and letting  $\epsilon \to 0$, we get
\begin{equation*}
\lim_{\epsilon \to 0}\frac{d}{dx}r(x_0)=\frac{u_1-u_0}{h}.
\end{equation*}
Note that, again the method reduces to forward difference formula of $u^{\prime}_x$ at $x=x_0$.

\subsubsection{N=2}
Consider the interpolant for $N=2$, we have
 \begin{equation}\label{eq11}
r(x)=\displaystyle\lambda_0\phi_0(x)+\lambda_1\phi_1(x)+\lambda_2\phi_2(x).
\end{equation}
Using the interpolation condition $r(x_k)=u_k, k= 0,1,2$ the interpolation matrix becomes a symmetric matrix with all diagonal entries  $1$.  Thus, we have
\begin{equation*}
\displaystyle\begin{pmatrix}
1&\displaystyle\frac{1}{1+\epsilon^2h^2}&\displaystyle\frac{1}{1+4\epsilon^2h^2}\\\displaystyle\frac{1}{1+\epsilon^2h^2}&1&\displaystyle\frac{1}{1+\epsilon^2h^2}\\\displaystyle\frac{1}{1+4\epsilon^2h^2}&\displaystyle\frac{1}{1+\epsilon^2h^2}&1
\end{pmatrix}
\begin{pmatrix}
\displaystyle\lambda_0\\\lambda_1\\\lambda_2
\end{pmatrix}=\displaystyle\begin{pmatrix}
u_0\\u_1\\u_2
\end{pmatrix},
\end{equation*}
where $h=x_{k+1}-x_k.$  Solving for $\lambda_k, k=0,1,2$, we get
\begin{eqnarray*}
\begin{aligned}
& \displaystyle\lambda_0=\frac{\big(1+4\epsilon^2h^2\big)\big(\big(2+9\epsilon^2h^2+4\epsilon^4h^4\big)u_0-4\big(1+\epsilon^2h^2\big)u_1+\big(2-\epsilon^2h^2\big)u_2\big)}{8\epsilon^4h^4\big(5+2\epsilon^2h^2\big)},\\
& \displaystyle\lambda_1=-\frac{\big(1+\epsilon^2h^2\big)\big(\big(1+4\epsilon^2h^2\big)u_0-2\big(1+3\epsilon^2h^2+2\epsilon^4h^4\big)u_1+\big(1+4\epsilon^2h^2\big)u_2\big)}{2\epsilon^4h^4\big(5+2\epsilon^2h^2\big)},\\
& \displaystyle\lambda_2=\frac{\big(1+4\epsilon^2h^2\big)\big(\big(2-\epsilon^2h^2\big)u_0-4\big(1+\epsilon^2h^2\big)u_1+\big(2+9\epsilon^2h^2+4\epsilon^4h^4\big)u_2\big)}{8\epsilon^4h^4\big(5+2\epsilon^2h^2\big)}.
\end{aligned}
\end{eqnarray*}
Differentiating the interpolant $r(x)$  with respect to $x$, we get
\begin{equation*}\label{eq12}
\frac{d}{dx}r(x)=-\frac{2\lambda_0(x-x_0)\epsilon^2}{{(1+\epsilon^2(x-x_0)^2})^{2}}-\frac{2\lambda_1(x-x_1)\epsilon^2}{{(1+\epsilon^2(x-x_1)^2})^{2}}-\frac{2\lambda_2(x-x_2)\epsilon^2}{{(1+\epsilon^2(x-x_2)^2})^{2}},
\end{equation*}
and evaluating at $x=x_1$, we have
\begin{equation}\label{eq13}
\frac{d}{dx} r(x_1)=\frac{2\epsilon^2h}{{(1+\epsilon^2h^2})^2}(\lambda_2-\lambda_0).\end{equation}
Substituting value of $\lambda_0,\lambda_2$ in \eqref{eq13}, we get
\begin{equation}\label{eqiq2}
\frac{d}{dx} r(x_1)=\dfrac{\big(1+4\epsilon^2h^2\big)}{2h\big(1+\epsilon^2h^2\big)^{2}}\big(u_2-u_0\big),
\end{equation}
and letting  $\epsilon \to 0$, we get
\begin{equation*}
\lim_{\epsilon \to 0}\frac{d}{dx}r(x_1)=\frac{u_2-u_0}{2h}.
\end{equation*}
Note that, again the method reduces to centred difference formula of $u^{\prime}_x$ at $x=x_1$. In this way, we can further construct the finite-difference approximation of higher derivatives using IQ- and IMQ-RBFs.

\section{Adaptive  RBF  Method for IVPs}
We consider the initial value problem of the form
\begin{equation} \label{Eq31}
\dfrac{du}{dt}=f(t,u), \,\, a \leq t \leq b,
\end{equation}
with initial condition 
\begin{equation}
u(a)=u_0,
\end{equation}
where we assume $u(t) \in C^{\infty}[a,b]$ and $f(t,u)$ is a class of $C^{\infty}$ function. We divide the interval $[a,b]$ in uniform way $t_n=a+nh,$ $n=0,1,...,N,$ where $h=(b-a)/N$ is the grid size. In this section, we derive the adaptive inverse multi-quadratic, inverse quadratic RBF methods such as Adam-Bashforth and Adam-Moulton methods, and their modifications.

\subsection{IMQ-RBF methods: Derivation and consistency analysis}
\subsubsection{IMQ-RBF AB2 method:}
We now derive two-step Adam-Bashforth method using adaptive IMQ-RBFs. For this, we approximate $f(t,u)$ using RBF interpolant $r(t)$ as
\begin{equation}
r(t_{n+i})=f(t_{n+i},u_{n+i})=f_{n+i},\,\,\, \text{for}: i=0,1.
\end{equation}
From equation \eqref{eq2},
\begin{equation}
r(t)=\frac{\lambda_0}{\sqrt{1+\epsilon_n^2(t-t_n)^2}}+\frac{\lambda_1}{\sqrt{1+\epsilon_n^2(t-t_{n+1})^2}},
\end{equation}
where,
\begin{equation*}
\begin{split}
&\lambda_0=\displaystyle\frac{(1+\epsilon_n^2h^2)}{\epsilon_n^2h^2}\bigg(f_n-\frac{f_{n+1}}{\sqrt{1+\epsilon_n^2h^2}}\bigg),\\
&\lambda_1=\displaystyle\frac{(1+\epsilon_n^2h^2)}{\epsilon_n^2h^2}\bigg(f_{n+1}-\frac{f_n}{\sqrt{1+\epsilon_n^2h^2}}\bigg).
\end{split}
\end{equation*}
Here, $r(t)$ is local approximation to right-hand-side (RHS) of equation \eqref{Eq31}. Considering the equation $v^{\prime}=r(t)$. Since, we have
\begin{equation*}
v(t_{n+2})-v(t_{n+1})=\int_{t_{n+1}}^{t_{n+2}} v^{\prime} \,dt \:= \: \int_{t_{n+1}}^{t_{n+2}} r(t) \,dt 
\end{equation*}
On simplifying, we get
\begin{equation}\label{IMQE-RBF}
v_{n+2}=v_{n+1}+\int_{t_{n+1}}^{t_{n+2}} r(t) \,dt 
\end{equation}
where,
\begin{equation*}
\begin{aligned}
\int_{t_{n+1}}^{t_{n+2}} r(t) \,dt= \lambda_0\int_{t_{n+1}}^{t_{n+2}}\displaystyle\frac{1}{\sqrt{1+\epsilon_n^2(t-t_n)^2}} \,dt+\lambda_1 \int_{t_{n+1}}^{t_{n+2}}\displaystyle\frac{1}{\sqrt{1+\epsilon_n^2(t-t_{n+1})^2}} \,dt.
\end{aligned}
\end{equation*}
Now, we compute
\begin{equation*}
\begin{split}
\int_{t_{n+1}}^{t_{n+2}}\displaystyle\frac{1}{\sqrt{1+\epsilon_n^2(t-t_n)^2}} \,dt=& \displaystyle\frac{1}{\epsilon_n}\bigg(\sinh^{-1}\big(2\epsilon_nh\big)-\sinh^{-1}\big(\epsilon_nh\big)\bigg),\\
\int_{t_{n+1}}^{t_{n+2}}\displaystyle\frac{1}{\sqrt{1+\epsilon_n^2(t-t_{n+1})^2}} \,dt =& \displaystyle\frac{1}{\epsilon_n}\bigg(\sinh^{-1}\big(\epsilon_nh\big)\bigg).
\end{split}
\end{equation*}
From equation \eqref{IMQE-RBF},  the modified method of the form
\begin{equation}\label{eq316}
v_{n+2}=v_{n+1}+h\big(\beta_0f_n+\beta_1f_{n+1}\big)
\end{equation}
where,
\begin{eqnarray*}
\begin{aligned}
&\beta_0=\displaystyle\frac{\sqrt{1+\epsilon_n^2h^2}}{\epsilon_n^3h^3}\bigg(\sqrt{1+\epsilon_n^2h^2}\sinh^{-1}(2\epsilon_nh)-\big(1+\sqrt{1+\epsilon_n^2h^2}\big)\sinh^{-1}(\epsilon_nh)\bigg),\\
& \beta_1=\displaystyle\frac{\sqrt{1+\epsilon_n^2h^2}}{\epsilon_n^3h^3}\bigg(\big(1+\sqrt{1+\epsilon_n^2h^2}\big)\sinh^{-1}(\epsilon_nh)-\sinh^{-1}(2\epsilon_nh)\bigg).
\end{aligned}
\end{eqnarray*}
Using Taylor series expansion around point $t=t_n$ to find local truncation error, the modified Adam-Bashforth method is then given by  
\begin{equation}\label{trunc}
\begin{split}
\tau_n&=\displaystyle\frac{u_{n+2}-u_{n+1}}{h}-\big(\beta_0f_n+\beta_1f_{n+1}\big),\\
&=\bigg(\frac{5}{12}\epsilon_n^2u_n^\prime+\displaystyle\frac{5}{12}u_n^{(3)}\bigg)h^2+\bigg(\displaystyle\frac{41}{24}\epsilon_n^2u_n^{\prime\prime}+\displaystyle\frac{3}{8}u_n^{(4)}\bigg)h^3+O(h^4).
\end{split}
\end{equation}
Note that, the modified Adam-Bashforth-two-point method with the IMQ-RBF  yields the second order accuracy as the leading error term is of $O(h^2)$. 
However, the coefficient of the second order term is not uniquely determined as it contains the free shape parameter $\epsilon_n.$\\
\newline
For $O(h^3)$: If we allow the leading error term to be zero, then we can eliminate the first term in the local truncation error so that we arrive at  third order of convergence. Thus,
$$\bigg(\displaystyle\frac{5}{12}\epsilon_n^2u_n^\prime+\displaystyle\frac{5}{12}u_n^{(3)}\bigg)=0$$ yields,
\begin{equation}\label{second}
\epsilon_n^2=-\frac{u_n^{(3)}}{u_n'}.
\end{equation}
The solution varies with the value of $n$ due to the presence of index n in the value of optimal shape parameter, $\epsilon_n$.\\
\newline
 For $O(h^4)$: If we allow the two leading error terms to be zero in local truncation error, then we could arrive at  fourth order of convergence. Thus,
$$ \bigg(\displaystyle\frac{5}{12}\epsilon_n^2u_n^\prime+\displaystyle\frac{5}{12}u_n^{(3)}\bigg)h^2+\bigg(\displaystyle\frac{41}{24}\epsilon_n^2u_n^{\prime\prime}+\frac{3}{8}u_n^{(4)}\bigg)h^3=0,$$
yields
\begin{equation}
\epsilon_n^2=-\displaystyle \dfrac{10u_n^{(3)}+9hu_n^{(4)}}{10u_n^\prime+41hu_n^{\prime\prime}}.
\end{equation}
As step size $h \to 0$, we get
\begin{equation}
\lim_{h \to 0}\epsilon_n^2=\lim_{h \to 0}\epsilon_n^2=-\displaystyle \dfrac{10u_n^{(3)}
+9hu_n^{(4)}}{10u_n^\prime+41hu_n^{\prime\prime}}
\\
=-\frac{u_n^{(3)}}{u_n'},
\end{equation}
which is identical to the value of $\epsilon_n^2$ in \eqref{second}. For achieving fifth and higher-order accuracy a complicated algebraic procedure needs to be employed to determine the optimal shape parameter $\epsilon_n^2$.
Now, expanding the coefficients of $f_n$ and $f_{n+1}$ in \eqref{eq316} using Taylor series about point $t_n$, and  discarding the higher order terms, we get the modified RBF-IMQ-AB2 method as
\begin{equation}\label{eqAB2imq}
v_{n+2}=v_{n+1}+\bigg(-\displaystyle\frac{h}{2}+\frac{31}{24}\epsilon_n^2h^3\bigg)f_n+\bigg(\displaystyle\frac{3h}{2}-\frac{41}{24}\epsilon_n^2h^3\bigg)f_{n+1}.
\end{equation}
The corresponding local truncation error is
\begin{equation}
\begin{split}
\tau_n=& \displaystyle\frac{u_{n+2}-u_{n+1}}{h}-\bigg(-\displaystyle\frac{1}{2}+\frac{31}{24}\epsilon_n^2h^2\bigg)f_n-\bigg(\frac{3}{2}-\frac{41}{24}\epsilon_n^2h^2\bigg)f_{n+1}\\
=& \bigg(\frac{5}{12}\epsilon_n^2u_n^\prime+\frac{5}{12}u_n^{(3)}\bigg)h^2+\bigg(\frac{41}{24}\epsilon_n^2u_n^{\prime\prime}+\frac{3}{8}u_n^{(4)}\bigg)h^3+O(h^4)
\end{split}
\end{equation}
Thus, the optimal value of shape parameter is $\epsilon_n$=$-\displaystyle\frac{u_n^{(3)}}{u_n^\prime}$. If we replace the third derivative $u_n^{(3)}=f_n''$ with the centred difference formula $\displaystyle\frac{f_{n+1}-2f_n+f_{n-1}}{h^2}$, we get
\begin{equation}\label{1q}
\epsilon_n^2=-\dfrac{f_{n+1}-2f_n+f_{n-1}}{h^2f_n}.
\end{equation}
 The obtained $\epsilon_n^2$ in \eqref{1q} yields the same desired order of accuracy. Now, we substitute this in equation \eqref{eqAB2imq} and use the obtained result to compute accuracy and rate of convergence of numerical solutions of some IVPs, which are presented in numerical section.

\subsubsection{IMQ-RBF-AB3 method:}
We can also derive IMQ-RBF three-step Adam-Bashforth method as
\begin{equation}\label{eqAB3imq}
v_{n+3}=v_{n+2}+\bigg(\displaystyle\frac{5h}{12}-\frac{213}{80}\epsilon_n^2h^3\bigg)f_n+\bigg(-\displaystyle\frac{4h}{3}+\frac{87}{10}\epsilon_n^2h^3\bigg)f_{n+1}+\bigg(\displaystyle\frac{23h}{12}-\frac{483}{80}\epsilon_n^2h^3\bigg)f_{n+2}.
\end{equation}
Considering the local truncation analysis, the optimal shape parameter is given by
\begin{equation}\label{second3}
\epsilon_n^2=-\frac{u_n^{(4)}}{9u_n^{\prime\prime}}.
\end{equation}
If we replace the fourth derivative $u_n^{(4)}=f_n'''$ and second derivative $u_n''=f_n'$ with the backward difference formula $\displaystyle\frac{f_{n+2}-3f_{n+1}+3f_n-f_{n-1}}{h^3}$ and $\displaystyle\frac{f_n-f_{n-1}}{h}$ respectively, we get
\begin{equation}
\epsilon_n^2=-\dfrac{f_{n+2}-3f_{n+1}+3f_{n}-f_{n-1}}{9h^2(f_n-f_{n-1})}.
\end{equation}

\subsubsection{IMQ-RBF-AM2 method:}
The RBF one-step Adams-Moulton method might be obtained in a similar manner. Since
$$ v\left(t_{n+1}\right)-v\left(t_{n}\right)=\int_{t_{\mathrm{n}}}^{t_{n+1}} v^{\prime} \mathrm{d} t=\int_{t_{n}}^{t_{n+1}} r(t) \mathrm{d}t,$$
we have
$$ v_{n+1}=v_{n}+\int_{t_{n}}^{t_{n+1}} r(t) dt,$$
and
$$ \int_{t_{n}}^{t_{n+1}} r(t) d t=\lambda_{0} \int_{t_{n}}^{t_{n+1}} \frac{1}{\sqrt{1+\epsilon_{n}^{2}\left(t-t_{n}\right)^{2}}} d t+\lambda_{1} \int_{t_{n}}^{t_{n+1}} \frac{1}{\sqrt{1+\epsilon_{n}^{2}\left(t-t_{n+1}\right)^{2}}} dt,$$
After expanding integrals, we obtain
$$
\begin{aligned}
\int_{t_{n}}^{t_{n+1}} \frac{1}{\sqrt{1+\epsilon_{n}^{2}\left(t-t_{n}\right)^{2}}} d t &=\frac{\sinh ^{-1}(\epsilon_{n} h)}{\epsilon_{n} },
\end{aligned}
$$
and
$$
\begin{aligned}
\int_{t_{n}}^{t_{n+1}} \frac{1}{\sqrt{1+\epsilon_{n}^{2}\left(t-t_{n+1}\right)^{2}}} d t &=\frac{\sinh ^{-1}(\epsilon_{n} h)}{\epsilon_{n} }.
\end{aligned}
$$
As a result, the revised scheme is as follows: 
$$
v_{n+1}=v_{n}+h \beta\left(f_{n+1}+f_{n}\right),
$$
where
$$
\beta=\frac{\left(\epsilon_{n}^2 h^2-\sqrt{\epsilon_{n} ^2 h^2+1}+1\right) \sinh ^{-1}(\epsilon_{n}  h)}{\epsilon_{n} ^3 h^3}.
$$
Hence,
\begin{equation}\label{eqn:ma1}
v_{n+2}=v_{n+1}+h \beta\left(f_{n+2}+f_{n+1}\right)
\end{equation}
where
$$ \beta=\frac{\left(\epsilon_{n+1}^2 h^2-\sqrt{\epsilon_{n+1} ^2 h^2+1}+1\right) \sinh ^{-1}(\epsilon_{n+1}  h)}{\epsilon_{n+1} ^3 h^3}.$$
Using the Taylor Series expansion around $t=t_{n}$, the local truncation error of $\eqref{eqn:ma1}$ is given by 
$$
\begin{aligned}
\tau_{n+1} &=\frac{u_{n+2}-u_{n+1}}{h}-\beta\left(f_{n+2}+f_{n+1}\right) \\
&=\left(-\frac{1}{12} \epsilon_{n+1}^{2} u_{n}^{\prime}-\frac{1}{12} u_{n}^{(3)}\right) h^{2}+\left(-\frac{1}{8} \epsilon_{n+1}^{2} u_{n}^{\prime \prime}-\frac{1}{8} u_{n}^{(4)}\right) h^{3}+O\left(h^{4}\right) .
\end{aligned}
$$
The optimal value of $\epsilon_{n+1}^2$ to attain third order accuracy is then calculated as follows: 
\begin{equation}\label{eqn:Ma3}
\epsilon_{n+1}^{2}=-\frac{u_{n}^{(3)}}{u_{n}^{\prime}},
\end{equation}
The optimum value of $\epsilon_{n+1}^2$ for fourth order accuracy is provided by the following equation. 
\begin{equation*}
\epsilon_{n+1}^{2}=-\frac{2 u_{n}^{(3)}+3 h u_{n}^{(4)}}{2 u_{n}^{\prime}+3 h u_{n}^{\prime \prime}}.
\end{equation*}
For the approximation form of $\epsilon_{n+1}^{2}$, we further extend $\beta$ in the Taylor series, yielding 
\begin{equation*}
v_{n+2}=v_{n+1}+\left(\frac{h}{2}+\frac{\epsilon_{n+1}^{2} h^{3}}{24}+O\left(h^{5}\right)\right)\left(f_{n+2}+f_{n+1}\right)
\end{equation*}
After we remove the higher order term from the above equation, we get 
\begin{equation}\label{eqAM2imq}
{v_{n+2}=v_{n+1}+\left(\frac{h}{2}+\frac{\epsilon_{n+1}^{2} h^{3}}{24}\right)\left(f_{n+2}+f_{n+1}\right)}.
\end{equation}
The corresponding local truncation error of $\eqref{eqAM2imq}$ is
$$
\begin{aligned}
\tau_{n+1} &=\frac{u_{n+2}-u_{n+1}}{h}-\left(\frac{1}{2}+\frac{\epsilon_{n+1}^{2} h^{2}}{24}\right)\left(f_{n+2}+f_{n+1}\right) \\
&=\left(-\frac{1}{12} \epsilon_{n+1}^{2} u_{n}^{\prime}-\frac{1}{12} u_{n}^{(3)}\right) h^{2}+\left(-\frac{1}{8} \epsilon_{n+1}^{2} u_{n}^{\prime \prime}-\frac{1}{8} u_{n}^{(4)}\right) h^{3}+O\left(h^{4}\right)
\end{aligned}
$$
Then the best value of $\epsilon_{n+1}^{2}$ for achieving the third order of accuracy is still $\eqref{eqn:Ma3}$.
The central difference $\dfrac{f_{n+1}-2 f_{n}+f_{n-1}}{h^{2}}$ replaces the third derivative $u_n^{(3)}=f_{n}^{\prime\prime}$ gives,
\begin{equation}
{\epsilon_{n+1}^{2}=-\frac{f_{n+1}-2 f_{n}+f_{n-1}}{h^{2} f_{n}}}
\end{equation}
Note that the resulting $\epsilon_{n+1}^{2}$ has the same accuracy of desired order. 

\subsubsection{IMQ-RBF-AM3 method:}
Similar to above derivations, the IMQ-RBF two-step Adam-Moulton(AM-2)  method is given by
\begin{equation}\label{eqAM3imq}
v_{n+3}=v_{n+2}+\left(-\frac{h}{12}-\frac{3 \epsilon_{n+1}^2 h^3}{80}\right)f_{n+1}+\left(\frac{2 h}{3}-\frac{3 \epsilon_{n+1}^2 h^3}{10}\right) f_{n+2}+\left(\frac{5 h}{12}+\frac{27 \epsilon_{n+1}^2 h^3}{80}\right) f_{n+3},
\end{equation}
with optimal shape parameter is 
\begin{equation}\label{eqn:Ma3}
\epsilon_{n+1}^{2}=-\frac{u_{n}^{(4)}}{9u_{n}^{(2)}}.
\end{equation}
Further, we replace the fourth derivative $u_{n}^{(4)}=f_{n}^{(3)}$ and the second derivative $u_{n}^{(2)}=f_{n}^{\prime}$ with the finite difference formulas as
\begin{equation}
{\epsilon_{n+1}^{2}=-\frac{f_{n+2}-3 f_{n+1}+3 f_{n}-f_{n-1}}{9h^{2}(f_{n}-f_{n-1})}}.
\end{equation}

\subsection{IQ-RBF methods: Derivation and consistency analysis}
\subsubsection{IQ-RBF-AB2 method}
Now, we derive the two-step Adam-Bashforth method with inverse-quadratic RBFs. In a similar way to above derivation, we approximate $f(t,u)$ using RBF interpolant $r(t)$
\begin{equation}
r(t_{n+i})=f(t_{n+i},u_{n+i})=f_{n+i},\,\,\,\, \text{for} i=0,1
\end{equation}
From equation \eqref{eq5}, we have
\begin{equation}
r(t)=\frac{\lambda_0}{1+\epsilon_n^2(t-t_n)^2}+\frac{\lambda_1}{1+\epsilon_n^2(t-t_{n+1})^2},
\end{equation}
where,
\begin{equation*}
\begin{split}
&\lambda_0=\displaystyle\frac{(1+\epsilon_n^2h^2)}{(\epsilon_n^2h^2)(2+\epsilon^2h^2)}\bigg((1+\epsilon^2h^2)f_n-f_{n+1}\bigg),\\
&\lambda_1=\displaystyle\frac{(1+\epsilon_n^2h^2)}{(\epsilon_n^2h^2)(2+\epsilon^2h^2)}\big((1+\epsilon^2h^2)f_{n+1}-f_{n}\big).
\end{split}
\end{equation*}
Note that, $r(t)$ is an approximation to RHS of equation \eqref{Eq31}. Considering the equation $v^{\prime}=r(t)$, and from the fundamental theorem of Integral calculus, we have the following;
\begin{equation*}
v(t_{n+2})-v(t_{n+1})=\int_{t_{n+1}}^{t_{n+2}} v^{\prime} \,dt \:= \: \int_{t_{n+1}}^{t_{n+2}} r(t) \,dt. 
\end{equation*}
On simplifying, we get
\begin{equation}\label{IQE-RBF}
v_{n+2}=v_{n+1}+\int_{t_{n+1}}^{t_{n+2}} r(t) \,dt, 
\end{equation}
with,
\begin{equation*}
\int_{t_{n+1}}^{t_{n+2}} r(t) \,dt= \lambda_0\int_{t_{n+1}}^{t_{n+2}}\frac{1}{1+\epsilon_n^2(t-t_n)^2} \,dt+\lambda_1 \int_{t_{n+1}}^{t_{n+2}}\frac{1}{1+\epsilon_n^2(t-t_{n+1})^2} \,dt.
\end{equation*}
Now, we have
\begin{equation*}
\begin{split}
&\int_{t_{n+1}}^{t_{n+2}}\frac{1}{1+\epsilon_n^2(t-t_n)^2} \,dt =\frac{1}{\epsilon_n}\bigg(\tan^{-1}\big(2\epsilon_nh\big)-\tan^{-1}\big(\epsilon_nh\big)\bigg),\\
& \int_{t_{n+1}}^{t_{n+2}}\frac{1}{1+\epsilon_n^2(t-t_{n+1})^2} \,dt =\frac{1}{\epsilon_n}\bigg(\tan^{-1}\big(\epsilon_nh\big)\bigg).
\end{split}
\end{equation*}
Thus, the equation \eqref{IQE-RBF} turns in the form 
\begin{equation}\label{eq326}
v_{n+2}=v_{n+1}+h\big(\beta_0f_n+\beta_1f_{n+1}\big)
\end{equation}
where,
\begin{eqnarray*}
\begin{split}
\beta_0=& \displaystyle\frac{1+\epsilon_n^2h^2}{\epsilon_n^3h^3(2+\epsilon_n^2h^2)}\bigg(\big(1+\epsilon_n^2h^2\big)\tan^{-1}(2\epsilon_nh)-\big(2+\epsilon_n^2h^2\big)\tan^{-1}(\epsilon_nh)\bigg),\\
\beta_1=& \displaystyle\frac{1+\epsilon_n^2h^2}{\epsilon_n^3h^3(2+\epsilon_n^2h^2)}\bigg(\big(2+\epsilon_n^2h^2\big)\tan^{-1}(\epsilon_nh)-\tan^{-1}(2\epsilon_nh)\bigg).
\end{split}
\end{eqnarray*}
With the use of Taylor series expansion about point $t=t_n$ to find local truncation error, we get 
\begin{equation}\label{trunc2}
\begin{split}
\tau_n&=\frac{u_{n+2}-u_{n+1}}{h}-\big(\beta_0f_n+\beta_1f_{n+1}\big)\\
&=\bigg(\frac{5}{6}\epsilon_n^2u_n^\prime+\frac{5}{12}u_n^{(3)}\bigg)h^2+\bigg(\frac{29}{12}\epsilon_n^2u_n^{\prime\prime}+\frac{3}{8}u_n^{(4)}\bigg)h^3+O(h^4).
\end{split}
\end{equation}
Note that, the Adam-Bashforth two-point method with the IQ-RBF yields the second order accuracy as the leading error term is of $O(h^2)$. 
However, the coefficient of the second order term is not uniquely determined as it contains the free shape parameter $\epsilon_n.$\\
\newline
For $O(h^3)$: If we allow the leading error term to be zero, then we can eliminate the first term in the truncation error so that we arrive at  third order of convergence. Thus,
$\bigg(\frac{5}{6}\epsilon_n^2u_n^\prime+\frac{5}{12}u_n^{(3)}\bigg)=0$ yields,
\begin{equation}\label{second2}
\epsilon_n^2=-\frac{u_n^{(3)}}{2u_n'}.
\end{equation}
\newline
 For $O(h^4)$: If we allow the two leading error terms to be zero, then we could arrive at  fourth order of convergence. Thus, 
$$ \bigg(\frac{5}{6}\epsilon_n^2u_n^\prime+\frac{5}{12}u_n^{(3)}\bigg)h^2+\bigg(\frac{29}{12}\epsilon_n^2u_n^{\prime\prime}+\frac{3}{8}u_n^{(4)}\bigg)h^3=0,$$
yields
\begin{equation}
\epsilon_n^2=-\displaystyle \dfrac{10u_n^{(3)}
+9hu_n^{(4)}}{20u_n^\prime+58hu_n^{\prime\prime}}.
\end{equation}
As step size $h \to 0$, we get
\begin{equation}
\lim_{h \to 0}\epsilon_n^2=\lim_{h \to 0}\epsilon_n^2=-\displaystyle \dfrac{10u_n^{(3)}
+9hu_n^{(4)}}{20u_n^\prime+58hu_n^{\prime\prime}}
\\
=-\frac{u_n^{(3)}}{2u_n'},
\end{equation}
which is identical to the value of $\epsilon_n^2$ in \eqref{second2}. For achieving fifth and higher-order accuracy a complicated algebraic procedure needs to be employed to determine the optimal shape parameter $\epsilon_n^2$. 
Expanding the coefficients of $f_n$ and $f_{n+1}$ in \eqref{eq326} using Taylor series about point $t_n$, and  discarding the higher order terms, we get the RBF-IQ-AB2 method as
\begin{equation}\label{eqAB2iq}
v_{n+2}=v_{n+1}+\bigg(-\frac{h}{2}+\frac{19}{12}\epsilon_n^2h^3\bigg)f_n+\bigg(\frac{3h}{2}-\frac{29}{12}\epsilon_n^2h^3\bigg)f_{n+1}.
\end{equation}
The corresponding local truncation error is
\begin{equation}
\begin{split}
\tau_n=\frac{u_{n+2}-u_{n+1}}{h}-\bigg(-\frac{1}{2}+\frac{19}{12}\epsilon_n^2h^2\bigg)f_n-\bigg(\frac{3}{2}-\frac{29}{12}\epsilon_n^2h^2\bigg)f_{n+1}\\
=\bigg(\frac{5}{6}\epsilon_n^2u_n^\prime+\frac{5}{12}u_n^{(3)}\bigg)h^2+\bigg(\frac{29}{12}\epsilon_n^2u_n^{\prime\prime}+\frac{3}{8}u_n^{(4)}\bigg)h^3+O(h^4).
\end{split}
\end{equation}
Thus the optimal value of shape parameter $\epsilon_n$=$-\displaystyle\frac{u_n^{(3)}}{2u_n^\prime}.$ If we replace the third derivative $u_n^{(3)}=f_n''$ with the centred difference formula $\displaystyle\dfrac{f_{n+1}-2f_n+f_{n-1}}{h^2}$ we get,
\begin{equation}
\epsilon_n^2=-\dfrac{f_{n+1}-2f_n+f_{n-1}}{2h^2f_n}.
\end{equation}
 We substitute this in equation \eqref{eqAB2iq} to get the modified RBF IQ Adam-Bashforth-2-point method.
 
\subsubsection{IQ-RBF-AB3 method}
Now we derive 3-step Adam-Bashforth method.  We approximate $f(t,u)$ using RBF interpolant $r(t)$ as
\begin{equation}
r(t_{n+i})=f(t_{n+i},u_{n+i})=f_{n+i},\,\,\, \text{for}\,\,\, i=0,1.
\end{equation}
From equation \eqref{eq5},
\begin{equation}
r(t)=\frac{\lambda_0}{1+\epsilon_n^2(t-t_n)^2}+\frac{\lambda_1}{1+\epsilon_n^2(t-t_{n+1})^2}+\frac{\lambda_2 }{1+\epsilon_n^2(t-t_{n+2})^2},
\end{equation}
with
\begin{eqnarray*}
\begin{split}
& \displaystyle\lambda_0=\frac{\big(1+4\epsilon^2h^2\big)\big(\big(2+9\epsilon^2h^2+4\epsilon^4h^4\big)f_n-4\big(1+\epsilon^2h^2\big)f_{n+1}+\big(2-\epsilon^2h^2\big)f_{n+2}\big)}{8\epsilon^4h^4\big(5+2\epsilon^2h^2\big)},\\
& \displaystyle\lambda_1=-\frac{\big(1+\epsilon^2h^2\big)\big(\big(1+4\epsilon^2h^2\big)f_n-2\big(1+3\epsilon^2h^2+
2\epsilon^4h^4\big)f_{n+1}+\big(1+4\epsilon^2h^2\big)f_{n+2}\big)}{2\epsilon^4h^4\big(5+2\epsilon^2h^2\big)}.\\
& \displaystyle\lambda_2=\frac{\big(1+4\epsilon^2h^2\big)\big(\big(2-\epsilon^2h^2\big)f_n-4\big(1+\epsilon^2h^2\big)f_{n+1}+\big(2+9\epsilon^2h^2+4\epsilon^4h^4\big)f_{n+2}\big)}{8\epsilon^4h^4\big(5+2\epsilon^2h^2\big)}.
\end{split}
\end{eqnarray*}
Here, $r(t)$ is an approximation to RHS of equation \eqref{Eq31}. Considering the equation $v^{\prime}=r(t)$, we have
\begin{equation*}
v(t_{n+3})-v(t_{n+2})=\int_{t_{n+2}}^{t_{n+3}} v^{\prime} \,dt \:= \: \int_{t_{n+2}}^{t_{n+3}} r(t) \,dt.
\end{equation*}
On further simplifying, we get
\begin{equation}\label{IQE2-RBF}
v_{n+3}=v_{n+2}+\int_{t_{n+2}}^{t_{n+3}} r(t) \,dt, 
\end{equation}
with,
\begin{equation*}
\int_{t_{n+2}}^{t_{n+3}} r(t) \,dt= \lambda_0\int_{t_{n+2}}^{t_{n+3}}\displaystyle\frac{1}{1+\epsilon_n^2(t-t_n)^2} \,dt+\lambda_1 \int_{t_{n+2}}^{t_{n+3}}\displaystyle\frac{1}{1+\epsilon_n^2(t-t_{n+1})^2} \,dt+\lambda_2 \int_{t_{n+2}}^{t_{n+3}}\displaystyle\frac{1}{1+\epsilon_n^2(t-t_{n+2})^2} \,dt.
\end{equation*}
Here, we compute the following.
\begin{equation*}
\begin{split}
& \int_{t_{n+2}}^{t_{n+3}}\displaystyle\frac{1}{1+\epsilon_n^2(t-t_n)^2} \,dt=\displaystyle\frac{1}{\epsilon_n}\bigg(\tan^{-1}\big(3\epsilon_nh\big)-\tan^{-1}\big(2\epsilon_nh\big)\bigg),\\
& \int_{t_{n+2}}^{t_{n+3}}\displaystyle\frac{1}{1+\epsilon_n^2(t-t_{n+1})^2} \,dt=\displaystyle\frac{1}{\epsilon_n}\bigg(\tan^{-1}\big(2\epsilon_nh\big)-\tan^{-1}\big(\epsilon_nh\big)\bigg),\\
& \int_{t_{n+2}}^{t_{n+3}}\displaystyle\frac{1}{\sqrt{1+\epsilon_n^2(t-t_{n+2})^2}} \,dt=\displaystyle\frac{1}{\epsilon_n}\bigg(\tan^{-1}\big(\epsilon_nh\big)\bigg).
\end{split}
\end{equation*}
With these computations, the equation \eqref{IQE2-RBF} become in the form 
\begin{equation}\label{eq346}
v_{n+3}=v_{n+2}+h\big(\beta_0f_n+\beta_1f_{n+1}+\beta_2f_{n+2}\big).
\end{equation}
Now, using Taylor series expansion about point $t=t_n$ to find local truncation error, we get 
\begin{equation}\label{trunc4}
\begin{split}
\tau_n&=\displaystyle\frac{u_{n+3}-u_{n+2}}{h}-\big(\beta_0f_n+\beta_1f_{n+1}+\beta_2f_{n+2}\big)\\
&=\bigg(\frac{36}{8}\epsilon_n^2u_n^{\prime\prime}+\displaystyle\frac{3}{8}u_n^{(4)}\bigg)h^3+\bigg(\displaystyle\frac{174}{25}\epsilon_n^4u_n^{\prime}+\displaystyle\frac{1487}{150}\epsilon_n^2u_n^{(3)}+\frac{193}{360}u_n^{(5)}\bigg)h^4+O(h^5)
\end{split}
\end{equation}
The modified Adam-Bashforth-3-point method with the IQ-RBF  yields the third order accuracy as the leading error term is of $O(h^3)$.\\
\newline
For $O(h^4)$: If we allow the leading error term to be zero, then we can eliminate the first term in the local truncation error so that we arrive at fourth order of convergence. Thus,
$\bigg(\displaystyle\frac{36}{8}\epsilon_n^2u_n^{\prime\prime}+\displaystyle\frac{3}{8}u_n^{(4)}\bigg)=0$ yields,
\begin{equation}\label{second4}
\epsilon_n^2=-\frac{u_n^{(4)}}{12u_n^{\prime\prime}}.
\end{equation}
The solution varies with the value of $n$ due to the presence of index n in the value of Optimal Shape Parameter, $\epsilon_n$. 
\par
For achieving fifth and higher-order accuracy a complicated algebraic procedure needs to be employed to determine the optimal shape parameter $\epsilon_n^2$.  Expanding the coefficients of $f_n$,$f_{n+1}$ and $f_{n+2}$ in \eqref{eq346} using Taylor series about point $t_n$, and  discarding the higher order terms, we get the modified RBF IQ AB3 method as 
\begin{equation}\label{eqAB3iq}
v_{n+3}=v_{n+2}+\bigg(\displaystyle\frac{5h}{12}-\frac{949}{300}\epsilon_n^2h^3\bigg)f_n+\bigg(\displaystyle-\frac{4h}{3}+\frac{812}{75}\epsilon_n^2h^3\bigg)f_{n+1}+\bigg(\displaystyle\frac{23h}{12}-\frac{2299}{300}\epsilon_n^2h^3\bigg)f_{n+2}.
\end{equation}
The corresponding local truncation error is
\begin{equation}
\begin{split}
\tau_n=\displaystyle\frac{u_{n+3}-u_{n+2}}{h}-\bigg(\displaystyle\frac{5}{12}-\frac{949}{300}\epsilon_n^2h^2\bigg)f_n-\bigg(-\frac{4}{3}+\frac{812}{75}\epsilon_n^2h^2\bigg)f_{n+1}-\bigg(\frac{23}{12}-\frac{2299}{300}\epsilon_n^2h^2\bigg)f_{n+2}\\
=\bigg(\frac{36}{8}\epsilon_n^2u_n^{\prime\prime}+\displaystyle\frac{3}{8}u_n^{(4)}\bigg)h^3+\bigg(\displaystyle\frac{174}{25}\epsilon_n^4u_n^{\prime}+\displaystyle\frac{1487}{150}\epsilon_n^2u_n^{(3)}+\frac{193}{360}u_n^{(5)}\bigg)h^4+O(h^5)
\end{split}
\end{equation}
Thus, the optimal value of shape parameter $\epsilon_n$=$-\displaystyle\frac{u_n^{(4)}}{12u_n^{\prime\prime}}$. If we replace the fourth derivative $u_n^{(4)}=f_n'''$ and second derivative $u_n''=f_n'$ with the backward difference formula $\displaystyle\frac{f_{n+2}-3f_{n+1}+3f_n-f_{n-1}}{h^3}$ and $\displaystyle\frac{f_n-f_{n-1}}{h} $respectively, we get,
\begin{equation}
\epsilon_n^2=-\dfrac{f_{n+2}-3f_{n+1}+3f_{n}-f_{n-1}}{12h^2(f_n-f_{n-1})}.
\end{equation}
We substitute this in equation \eqref{eqAB3iq} and therefore, we obtain a modified RBF-IQ Adam Bashforth-3-point method.

\subsubsection{IQ-RBF AM2 method:} In a similar manner of previous derivations, the RBF one-step Adams-Moulton method might be obtained.
\begin{equation}\label{eqAM2iq}
{v_{n+2}=v_{n+1}+\left(\frac{h}{2}+\frac{\epsilon_{n+1}^{2} h^{3}}{12}\right)\left(f_{n+2}+f_{n+1}\right)}.
\end{equation}
The corresponding local truncation error of $\eqref{eqAM2iq}$ is
$$
\begin{aligned}
\tau_{n+1} &=\frac{u_{n+2}-u_{n+1}}{h}-\left(\frac{1}{2}+\frac{\epsilon_{n+1}^{2} h^{2}}{12}\right)\left(f_{n+2}+f_{n+1}\right) \\
&=\left(-\frac{1}{6} \epsilon_{n+1}^{2} u_{n}^{\prime}-\frac{1}{12} u_{n}^{(3)}\right) h^{2}+\left(-\frac{1}{4} \epsilon_{n+1}^{2} u_{n}^{\prime \prime}-\frac{1}{8} u_{n}^{(4)}\right) h^{3}+O\left(h^{4}\right). 
\end{aligned}
$$
Then the best value of $\epsilon_{n+1}^{2}$ for achieving the third order of accuracy is  $\epsilon_{n+1}^{2}=\dfrac{u_n^{(3)}}{2u_n^{\prime}}.$
Further, the central difference $\dfrac{f_{n+1}-2 f_{n}+f_{n-1}}{h^{2}}$ replaces the third derivative $u_n^{(3)}=f_{n}^{\prime\prime}$, thus we have,
\begin{equation}
{\epsilon_{n+1}^{2}=-\dfrac{f_{n+1}-2 f_{n}+f_{n-1}}{2h^{2} f_{n}}},
\end{equation}
the resulting $\epsilon_{n+1}^{2}$ has the same accuracy of desired order. 

\subsubsection{IQ-RBF AM3 method:} The IQ-RBF AM3 method is
\begin{equation}\label{eqAM3iq}
v_{n+3}=v_{n+2}+\left(-\frac{h}{12}-\frac{19 \epsilon_{n+1}^2 h^3}{300}\right)f_{n+1}+\left(\frac{2 h}{3}-\frac{28 \epsilon_{n+1}^2 h^3}{75}\right) f_{n+2}+\left(\frac{5 h}{12}+\frac{131 \epsilon_{n+1}^2 h^3}{300}\right) f_{n+3}
\end{equation}
The corresponding local truncation error of $\eqref{eqAM3iq}$ is
$$
\begin{aligned}
\tau_{n+1} &=\left(-\frac{12}{24} \epsilon_{n+1}^{2} u_{n}^{(2)}-\frac{1}{24} u_{n}^{(4)}\right) h^{3}+\frac{1}{225} h^4 \left(-54 \epsilon_{n+1}^4 u_{n}^{\prime}-267 \epsilon_{n+1}^2 u_{n}^{(3)}-20 u_{n}^{(5)}\right)+O\left(h^{5}\right).
\end{aligned}
$$
One can obtain the best optimal shape parameter by making the first term vanishes and further, we can use finite difference formulas for optimal shape parameter as
\begin{equation}
{\epsilon_{n+1}^{2}=-\frac{f_{n+2}-3 f_{n+1}+3 f_{n}-f_{n-1}}{12h^{2}(f_{n}-f_{n-1})}},
\end{equation}
thus, we obtain modified IQ-RBF-AM3 method.

\subsection{Stability}
We have developed various adaptive RBF-IMQ and RBF-IQ techniques and studied the consistency analysis in the preceding section. In practise, for a non-linear method, we study the  absolute stability \cite{Leveque2007finite}. The focus of our attention in this section will be to evaluate and find the region of absolute stability for all proposed methods. Consider
\begin{equation}
{
v^{\prime}=\lambda v.
}
\label{1}
\end{equation}
Note that, the linear equation \eqref{1} is used to establish this qualitative analysis. We derive the stability polynomial $\pi(\zeta ; z)$ where $\zeta$ is the roots of the characteristic polynomial and in general, its coefficients depend on the value of z.
\begin{definition}\label{def}
The set of points $z$ in the complex plane for which the polynomial $\pi(\zeta ; z)$, called the stability polynomial, satisfies the following criteria is known as the region of absolute stability (or the stability area)\\
(i) If $\zeta_{j}$ are the roots of the polynomial $\pi(\zeta ; z)$ then $\left|\zeta_{j}\right| \leqslant 1$ for each root.\\
(ii) If $\zeta_{j}$ is a repeated root, then $\left|\zeta_{j}\right|<1$.
\end{definition}
We compute the stability polynomial and based on this, we derive the stability region for each of the proposed methods.\\
\newline
\textbf{ RBF-IMQ-AB2 method:} Recall the equation \eqref{eqAB2imq},
\begin{equation}\label{2}
{v_{n+2} = v_{n+1} + \left(-\frac{h}{2}+\frac{31\epsilon_{n}^{2}h^{3}}{24}\right)f_{n}+ \left(\frac{3h}{2}-\frac{41\epsilon_{n}^{2}h^{3}}{24}\right)f_{n+1}}
\end{equation}
with
\begin{equation}\label{3}
{\epsilon_{n}^{2} =-\frac{v_{n}^{(3)}}{v_{n}^{\prime}}}
\end{equation}\\
From \eqref{1}, we can obtain the values of $v_{n}^{(3)}$, $v_{n}^{\prime}$, $f_{n}$, $f_{n+1}$ as
\begin{equation}
{v_{n}^{\prime} =\lambda v_{n} ;\hspace{0.5cm} v_{n+1}^{\prime} =\lambda v_{n+1}}; \hspace{0.5cm} \:{v_{n}^{\prime\prime} =\lambda^{2} v_{n}}; \hspace{0.5cm}\:{v_{n}^{(3)} =\lambda^{3} v_{n}}.
\end{equation}
Substituting the above equations in \eqref{2} and \eqref{3}, we get
\begin{equation}\label{4a}
{v_{n+2} = v_{n+1} + \left(-\frac{h}{2}+\frac{31\epsilon_{n}^{2}h^{3}}{24}\right)\lambda v_{n}+ \left(\frac{3h}{2}-\frac{41\epsilon_{n}^{2}h^{3}}{24}\right)\lambda v_{n+1}}
\end{equation}
\begin{equation}\label{4}
{\epsilon_{n}^{2} =-\lambda^{2}}
\end{equation}
Substituting $\epsilon_{n}^{2}$ value from \eqref{4} into \eqref{4a}, we get
\begin{equation}
{v_{n+2} = v_{n+1} + \left(-\frac{\lambda h}{2}-\frac{31(\lambda h)^{3}}{24}\right)v_{n}+ \left(\frac{3 \lambda h}{2}+\frac{41 (\lambda h)^{3}}{24}\right)v_{n+1}}
\end{equation}
After  regrouping of terms by considering the $z=h\lambda$, we get
\begin{equation*}{
v_{n+2}=\left(1+\frac{3z}{2}+\frac{41z^{3}}{24}\right) v_{n+1}+\left(-\frac{z}{2}-\frac{31z^{3}}{24}\right)v_{n}
}\end{equation*}
Let $\dfrac{v_{n+2}}{v_{n+1}}=\zeta$ , $\dfrac{v_{n+1}}{v_{n}}=\zeta$ , $\dfrac{v_{n+2}}{v_{n}}=\zeta^{2}$.
Hence, the stability polynomial is
\begin{equation}{
\pi(\zeta ; z)=\zeta^{2}-\left(1+\frac{3z}{2}+\frac{41z^{3}}{24}\right)\zeta+\left(\frac{z}{2}+\frac{31z^{3}}{24}\right)
}\end{equation}
The roots of the above polynomial say, $\zeta_1$ and $\zeta_2$ and the stability region can be found by applying the definition \eqref{def}. Now, we present the stability polynomial of each method and thus, one can obtain the stability regions in a similar way  to the RBF-IMQ-AB2 method.\\
\newline
\textbf{RBF-IQ-AB2 method:} The stability polynomial is
\begin{equation}{
\pi(\zeta ; z)=\zeta^{2}-\left(1+\frac{3z}{2}+\frac{29z^{3}}{24}\right)\zeta+\left(\frac{z}{2}+\frac{19z^{3}}{24}\right).
}\end{equation}
\textbf{RBF-IMQ-AB3 method:} We obtain the stability polynomial as
\begin{equation}{
\pi(\zeta ; z)=\zeta^{3}-\left(1+\frac{23z}{12}+\frac{483z^{3}}{720}\right)\zeta^{2}+\left(\frac{4z}{3}+\frac{87z^{3}}{90}\right)\zeta-\left(\frac{5z}{12}+\frac{213z^{3}}{720}\right)\zeta.
}\end{equation}
\textbf{RBF-IQ-AB3 method:} The stability polynomial is
\begin{equation}{
\pi(\zeta ; z)=\zeta^{3}-\left(1+\frac{23z}{12}+\frac{2299z^{3}}{3600}\right)\zeta^{2}+\left(\frac{4z}{3}+\frac{812z^{3}}{900}\right)\zeta-\left(\frac{5z}{12}+\frac{949z^{3}}{3600}\right)\zeta.
}\end{equation}
\textbf{RBF-IMQ-AM2 method:} The stability polynomial is
\begin{equation}{
\pi(\zeta ; z)=\left(1-\frac{z}{2}+\frac{z^{3}}{24}\right)\zeta-\left(1+\frac{z}{2}-\frac{z^{3}}{24}\right)
}\end{equation}
\textbf{RBF-IQ-AM2 method:} The stability polynomial is
\begin{equation}{
\pi(\zeta ; z)=\left(1-\frac{z}{2}+\frac{z^{3}}{24}\right)\zeta-\left(1+\frac{z}{2}-\frac{z^{3}}{24}\right)
}\end{equation}
\textbf{RBF-IMQ-AM3 method:} The stability polynomial is
\begin{equation}{
\pi(\zeta ; z)=\left(1-\frac{5z}{12}+\frac{27z^{3}}{80}\right)\zeta^{2}-\left(1+\frac{2z}{3}+\frac{3z^{3}}{90}\right)\zeta-\left(-\frac{z}{12}+\frac{3z^{3}}{720}\right).
}\end{equation}
\textbf{RBF-IQ-AM3 method:} The obtained stability polynomial is
\begin{equation}{
\pi(\zeta ; z)=\left(1-\frac{5z}{12}+\frac{131z^{3}}{3600}\right)\zeta^{2}-\left(1+\frac{2z}{3}+\frac{28z^{3}}{900}\right)\zeta-\left(-\frac{z}{12}+\frac{19z^{3}}{3600}\right).
}\end{equation}
\begin{figure}
  \centering
  \begin{minipage}[b]{0.42\textwidth}
    \includegraphics[width=\textwidth]{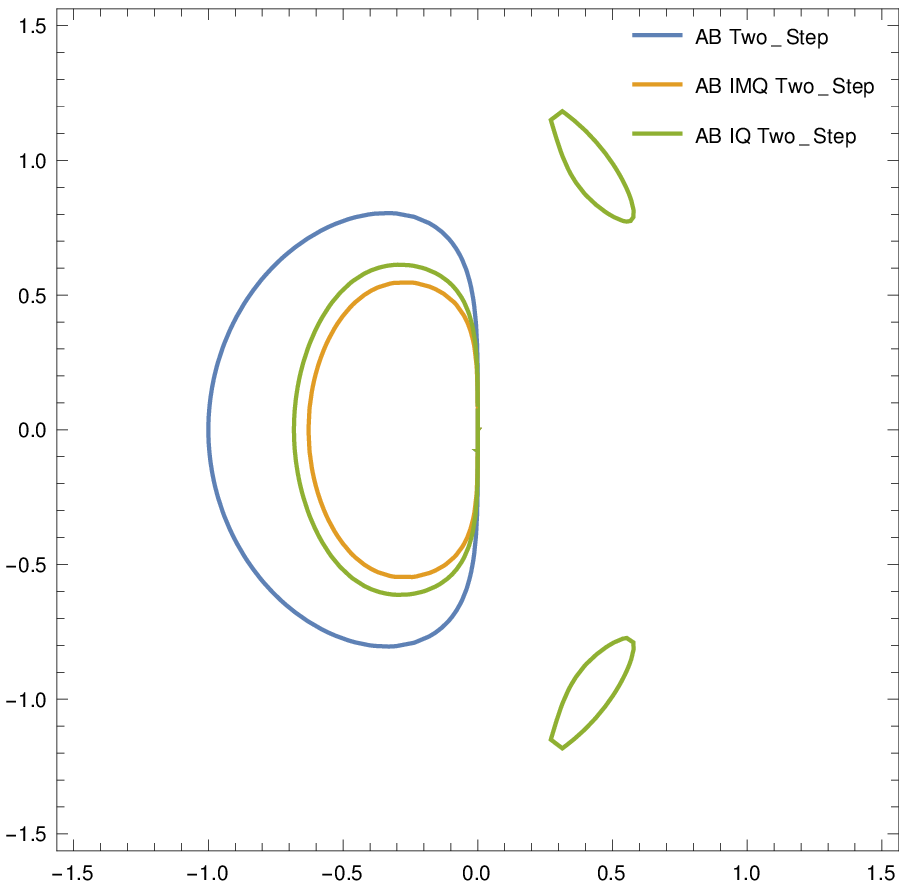}
    \caption{Stability region of all AB2 methods}\label{fig1}
  \end{minipage}
  \hfill
  \begin{minipage}[b]{0.42\textwidth}
    \includegraphics[width=\textwidth]{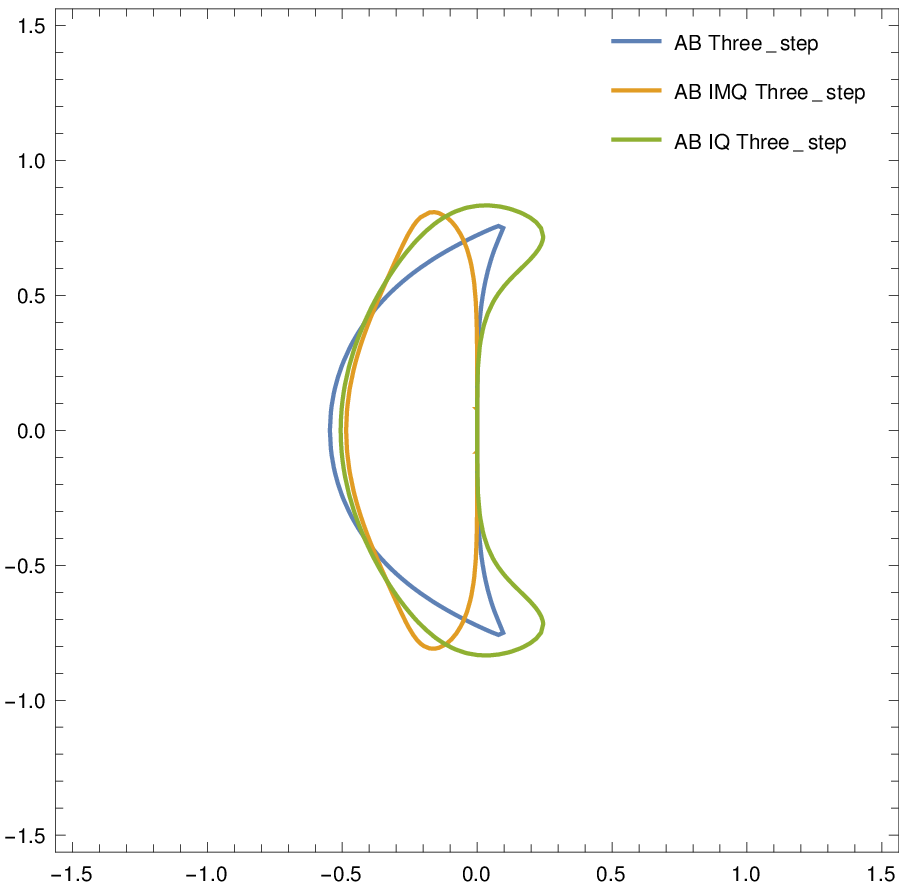}
    \caption{Stability region of all  AB3 methods}\label{fig2}
  \end{minipage}
\hfill
  \begin{minipage}[b]{0.42\textwidth}
    \includegraphics[width=\textwidth]{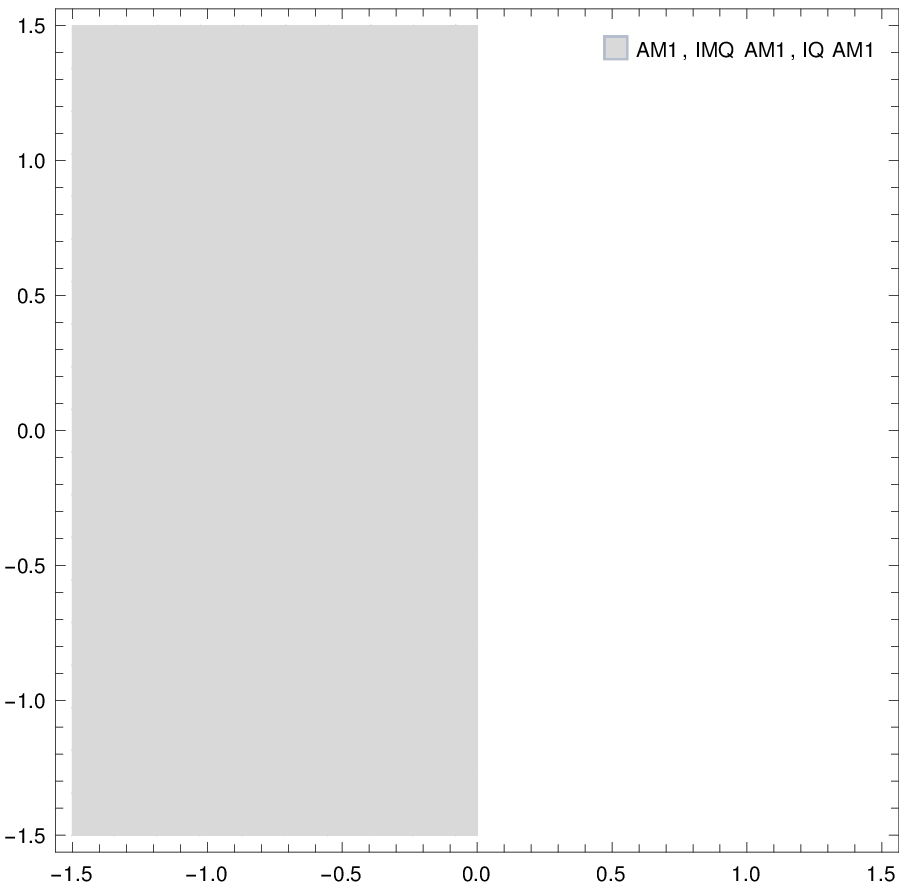}
    \caption{Stability region of all AM2 methods}\label{fig3}
  \end{minipage}
  \hfill
  \begin{minipage}[b]{0.42\textwidth}
     \includegraphics[width=\textwidth]{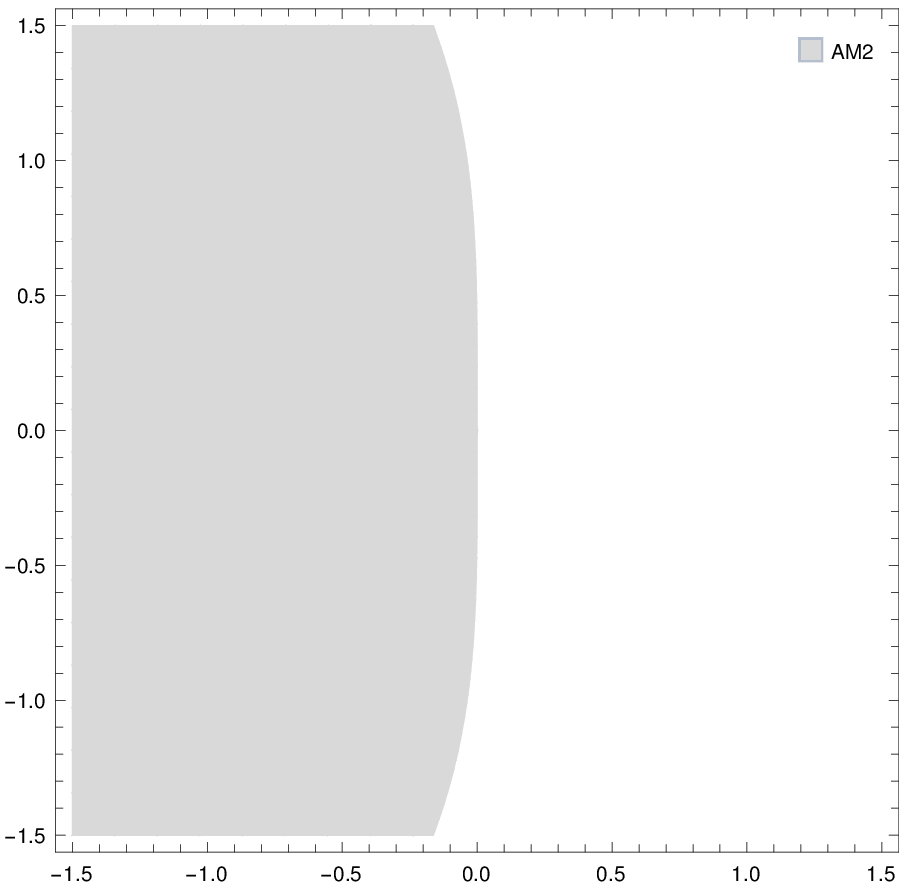}
    \caption{Stability region of AM3 method}\label{fig4}
  \end{minipage}
\hfill
  \begin{minipage}[b]{0.42\textwidth}
    \includegraphics[width=\textwidth]{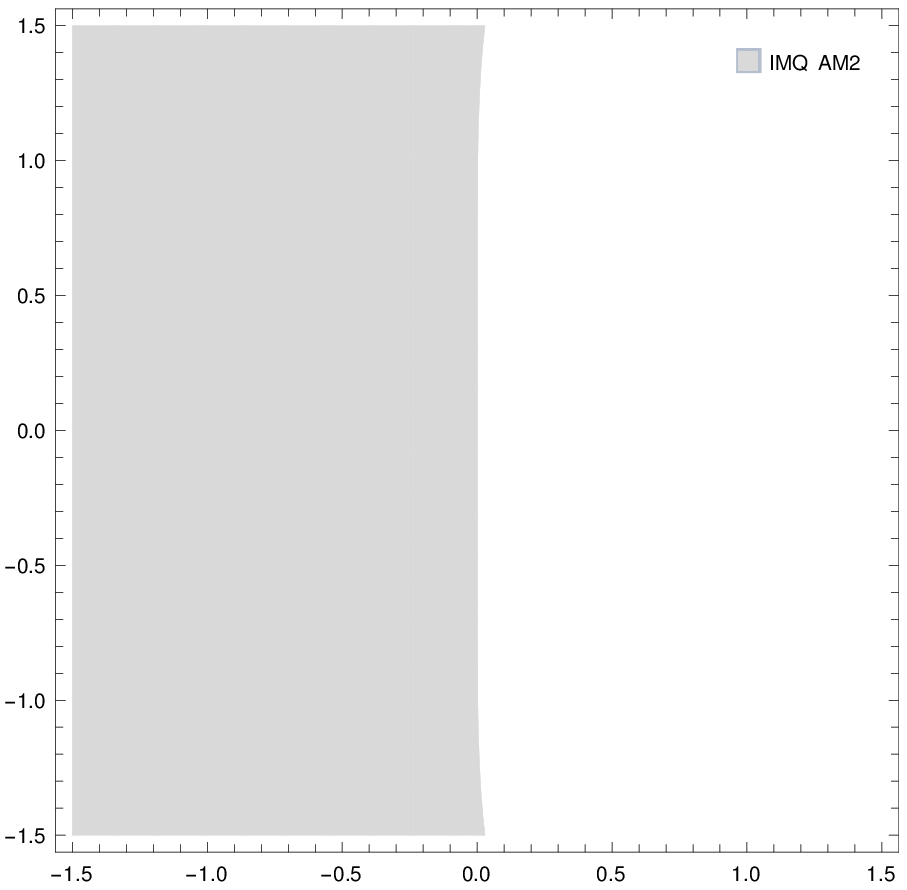}
    \caption{Stability region of RBF-IMQ -AM3 method}\label{fig5}
  \end{minipage}
    \hfill
  \begin{minipage}[b]{0.42\textwidth}
    \includegraphics[width=\textwidth]{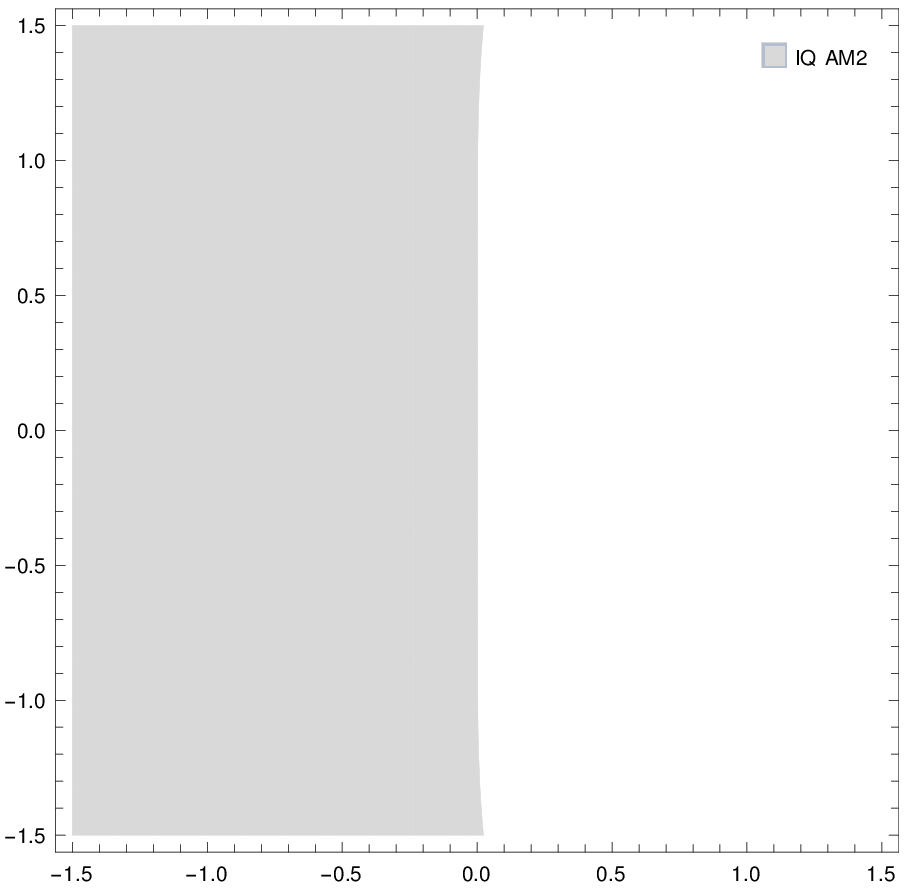}
    \caption{Stability region of RBF-IQ-AM3 method}\label{fig6}
  \end{minipage}
\end{figure}
\par
Figures \eqref{fig1} to \eqref{fig6} show the stability regions for various proposed methods. Figure \eqref{fig1} shows the original AB2 with proposed two-point RBF methods, and it is observed that the region of absolute stability for the proposed RBF methods is smaller than the original AB2 method. Thus, we see that although the proposed RBF methods yield the third order of accuracy compared to the second-order accuracy of the original AB2 method, they give better accuracy, but they do not perform better when stability is considered. The stability region shown in figure \eqref{fig2} contains the stability region of AB3, IMQ-RBF, and IQ-RBF methods. However, the proposed methods have an almost equivalent stability region, which is slightly bigger than the original AB3 point method; thus, the proposed methods perform better in terms of convergence and stability. Since the stability region of AM2 and proposed IMQ-AM2, IQ-AM2 shares the same regions, which is shown in figure \eqref{fig3}; thus, all methods perform well, but RBF methods yield a higher order of convergence.  In figures \eqref{fig4}, \eqref{fig5} and \eqref{fig6}, the stability regions of AM3, IMQ-AM3, IQ-AM3 methods and observed a similar behaviour as like AM2 methods. Thus, the modifications improvise the original method in terms of accuracy, while the stability remains unaffected.  In the following Table \ref{Table:Col}, we collect the AB2, AM3, RK2, RK3 and adaptive methods proposed here with order and optimal shape parameters. Number of function evaluations (FE), and the number of floating point operations (FPE) required to perform per step of the method are shown in Table \ref{Table:FE}.  Note that $\{n \pm, p \times, r \div\}$ denotes number of $n$-additions and$/$or subtractions, $p$-multiplications, $r$-divisions to be performed for a single iteration of the method, and the total number of floating point operations is calculated as $n+p+r$.  Note that to generate the numerical results, we have used the Python programming language on the machine having a 1.7 GHz Quad-Core Intel Core i5-4210U processor with 4GB memory.
\newpage
\begin{landscape}
\begin{table}[h!]
\begin{center}
\begin{tabular}{ l | l | l | p{4.5cm} }
\hline 
\textbf{Method} & \textbf{Numerical scheme} & \textbf{Order} & \textbf{Optimal $\epsilon_n^2$} \\
\hline
\textit{AB2}&$v_{n+2}=v_{n+1}+h\bigg(\dfrac{3}{2}f_{n+1}-\dfrac{1}{2}f_n\bigg)$&$O(h^2)$&\hspace{0.5cm}---\\
\textit{AM3}& $v_{n+2}=v_{n+1}+\dfrac{h}{2}\bigg(f_{n+2}+f_{n+1}\bigg)$ & $O(h^2)$&\hspace{0.5cm}---\\
\textit{RK2}&$v_{n+1}=v_{n}+\dfrac{(k_1+k_2)}{2},$&$O(h^2)$&\hspace{0.5cm}---\\ [0.5 ex]
  &$k_1=hf(v_n,t_n),  k_2=hf(v_n+k_1,t_n+h)$& & \\[0.5 ex]
  
 \textit{RK3}&$v_{n+1}=v_{n}+\dfrac{(k_1+4k_2+k_3)}{6},$&$O(h^3)$&\hspace{0.5cm}---\\ [0.5 ex]
  &$k_1=hf(v_n,t_n),  k_2=f(v_n+k_1/2,t_n+h/2), k_3=f(v_n+2k_2-k_1,t_n+h), $& & \\[0.5 ex] 
  
\hline  
\textit{IMQ-RBF AB2}  & $v_{n+2}=v_{n+1}+\bigg(-\displaystyle\frac{h}{2}+\frac{31}{24}\epsilon_n^2h^3\bigg)f_n+\bigg(\displaystyle\frac{3h}{2}-\frac{41}{24}\epsilon_n^2h^3\bigg)f_{n+1}$ & $O(h^3)$ & $-\dfrac{f_{n+1}-2f_n+f_{n-1}}{h^2f_n}$\\
\textit{IQ-RBF AB2}  & $v_{n+2}=v_{n+1}+\bigg(-\dfrac{h}{2}+\dfrac{19}{12}\epsilon_n^2h^3\bigg)f_n+\bigg(\dfrac{3h}{2}-\dfrac{29}{12}\epsilon_n^2h^3\bigg)f_{n+1}$ & $O(h^3)$ & $-\dfrac{f_{n+1}-2f_n+f_{n-1}}{2h^2f_n}$\\

 \textit{IMQ-RBF AM2}  & ${v_{n+2}=v_{n+1}+\left(\dfrac{h}{2}+\dfrac{\epsilon_{n+1}^{2} h^{3}}{24}\right)\left(f_{n+2}+f_{n+1}\right)}$ & $O(h^3) $& $-\dfrac{f_{n+1}-2 f_{n}+f_{n-1}}{h^{2} f_{n}}$\\
   \textit{IQ-RBF AM2}  & ${v_{n+2}=v_{n+1}+\left(\dfrac{h}{2}+\dfrac{\epsilon_{n+1}^{2} h^{3}}{12}\right)\left(f_{n+2}+f_{n+1}\right)}$ & $O(h^3) $& $-\dfrac{f_{n+1}-2 f_{n}+f_{n-1}}{2h^{2} f_{n}}$\\
 \hline    
  \textit{IMQ-RBF AB3}  & $v_{n+3}=v_{n+2}+\bigg(\displaystyle\frac{5h}{12}-\dfrac{213}{80}\epsilon_n^2h^3\bigg)f_n+\bigg(-\displaystyle\frac{4h}{3}+\dfrac{87}{10}\epsilon_n^2h^3\bigg)f_{n+1}+\bigg(\displaystyle\dfrac{23h}{12}-\dfrac{483}{80}\epsilon_n^2h^3\bigg)f_{n+2}$ & $O(h^4)$ & $-\dfrac{f_{n+2}-3f_{n+1}+3f_{n}-f_{n-1}}{9h^2(f_n-f_{n-1})}$\\
\textit{IQ-RBF AB3}  & $v_{n+3}=v_{n+2}+\bigg(\displaystyle\dfrac{5h}{12}-\frac{949}{300}\epsilon_n^2h^3\bigg)f_n+\bigg(\displaystyle-\dfrac{4h}{3}+\frac{812}{75}\epsilon_n^2h^3\bigg)f_{n+1}+\bigg(\displaystyle\dfrac{23h}{12}-\frac{2299}{300}\epsilon_n^2h^3\bigg)f_{n+2}
$ & $O(h^4)$ & $-\dfrac{f_{n+2}-3f_{n+1}+3f_{n}-f_{n-1}}{12h^2(f_n-f_{n-1})}$\\ 
 \textit{IMQ-RBF AM2}  & $v_{n+3}=v_{n+2}+\left(-\dfrac{h}{12}-\dfrac{3 \epsilon_{n+1}^2 h^3}{80}\right)f_{n+1}+\left(\dfrac{2 h}{3}-\dfrac{3 \epsilon_{n+1}^2 h^3}{10}\right) f_{n+2}+\left(\dfrac{5 h}{12}+\dfrac{27 \epsilon_{n+1}^2 h^3}{80}\right) f_{n+3}$ & $O(h^4) $& $-\dfrac{f_{n+2}-3 f_{n+1}+3 f_{n}-f_{n-1}}{9h^{2}(f_{n}-f_{n-1})}$\\
 \textit{IQ-RBF AM2}  & $v_{n+3}=v_{n+2}+\left(-\dfrac{h}{12}-\dfrac{19 \epsilon_{n+1}^2 h^3}{300}\right)f_{n+1}+\left(\dfrac{2 h}{3}-\dfrac{28 \epsilon_{n+1}^2 h^3}{75}\right) f_{n+2}+\left(\dfrac{5 h}{12}+\dfrac{131 \epsilon_{n+1}^2 h^3}{300}\right) f_{n+3}$ & $O(h^4) $& $-\dfrac{f_{n+2}-3 f_{n+1}+3 f_{n}-f_{n-1}}{12h^{2}(f_{n}-f_{n-1})}$\\
\hline
\end{tabular}
\caption{Summary of AB2, AM3, RK2, RK3, and proposed methods IQ-, and IMQ- AB2, AM2, AM3 Methods}
\label{Table:Col}
\end{center}
\end{table}
\end{landscape}

\begin{table}[h!]
\begin{center}
\begin{tabular}{ l | l | l p{3.5cm} }
\hline 
\textbf{Method} & \textbf{FE}  & \textbf{FPE} \\
\hline
 \textit{AB2}  &  2 &5=\{$2\pm, 2\times,1\div $\}\\
 \textit{AM2}&2&5=\{$2\pm, 2\times,1\div $\}\\
 \textit{RK2}&2&5=\{$2\pm, 2\times,1\div $\}\\
\textit{AB3}&3&8=\{$3\pm, 4\times,1\div $\}\\
\textit{AM3}&3&8=\{$3\pm, 4\times,1\div $\}\\
\textit{RK3}&3&6=\{$3\pm, 2\times,1\div $\}\\
\textit{IMQ-RBF AB2}  & 3&12=\{$4\pm, 6\times,2\div $\}\\
\textit{IQ-RBF AB2}  & 3&12=\{$4\pm,6\times,2\div $\}\\ 
\textit{IMQ-RBF AB3}  & 4&24=\{$6\pm, 12\times,6\div $\}\\
\textit{IQ-RBF AB3}  & 4&24=\{$6\pm, 12\times,6\div $\}\\ 
\textit{IMQ-RBF AM2}  & 4&12=\{$4\pm, 6\times,2\div $\}\\
\textit{IQ-RBF AM2}  & 4&12=\{$4\pm, 6\times,2\div $\}\\ 
\textit{IMQ-RBF AM3}  & 5&24=\{$6\pm, 10\times,6\div $\}\\
\textit{IQ-RBF AM3}  & 5&24=\{$6\pm, 10\times,6\div $\}\\ 
\hline
\end{tabular}
\caption{Number of function evaluations (FE)and floating point operations (FPE) to perform single iteration of method}
\label{Table:FE}
\end{center}
\end{table}

\section{Numerical results}
\begin{figure}[h!]
\centering
\subfloat{{\includegraphics[trim=1mm  1mm 0.2cm 0cm, clip=true, scale=0.5]{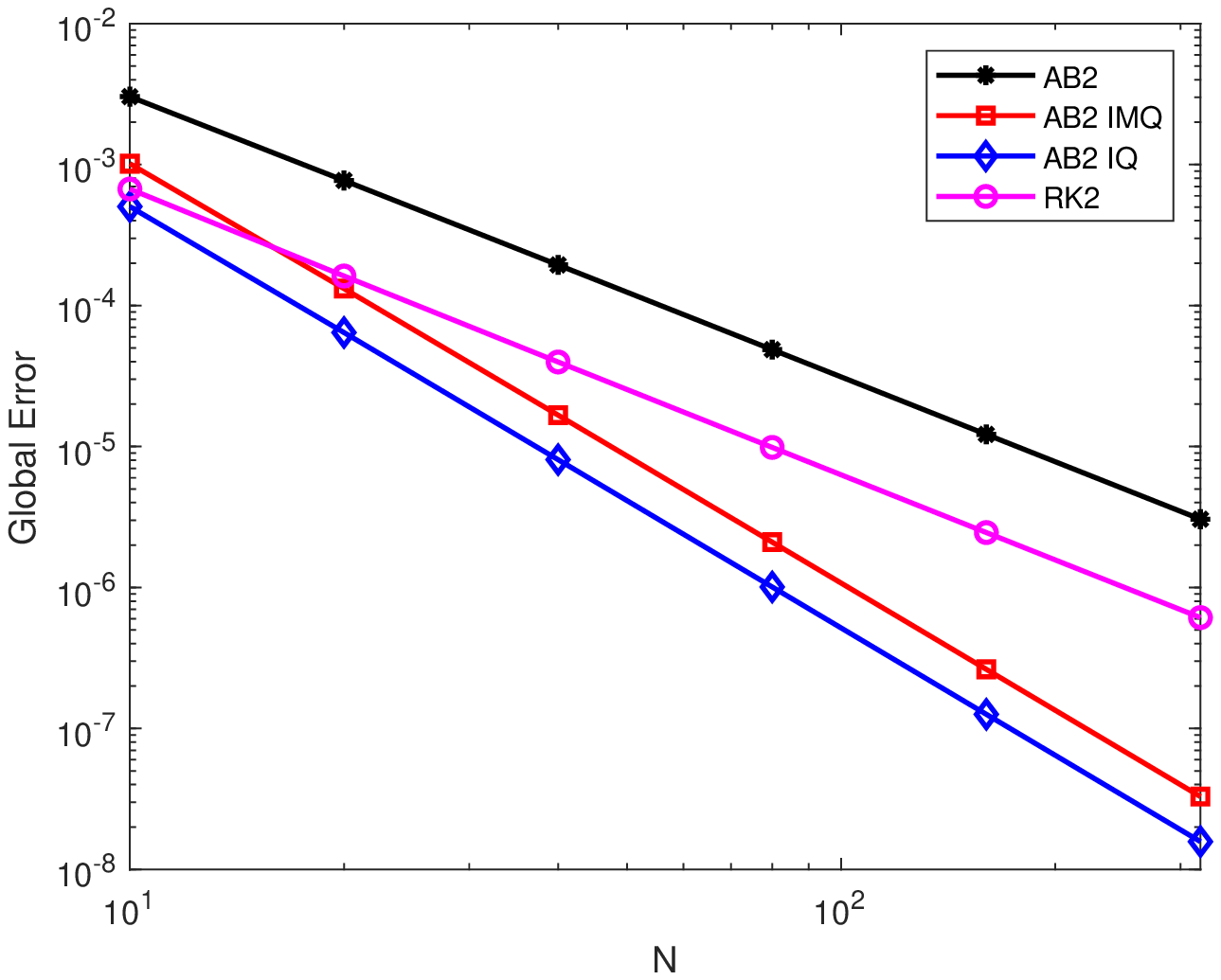} }}
\qquad
\subfloat{{\includegraphics[trim=1mm  1mm 0.2cm 0cm, clip=true, scale=0.5]{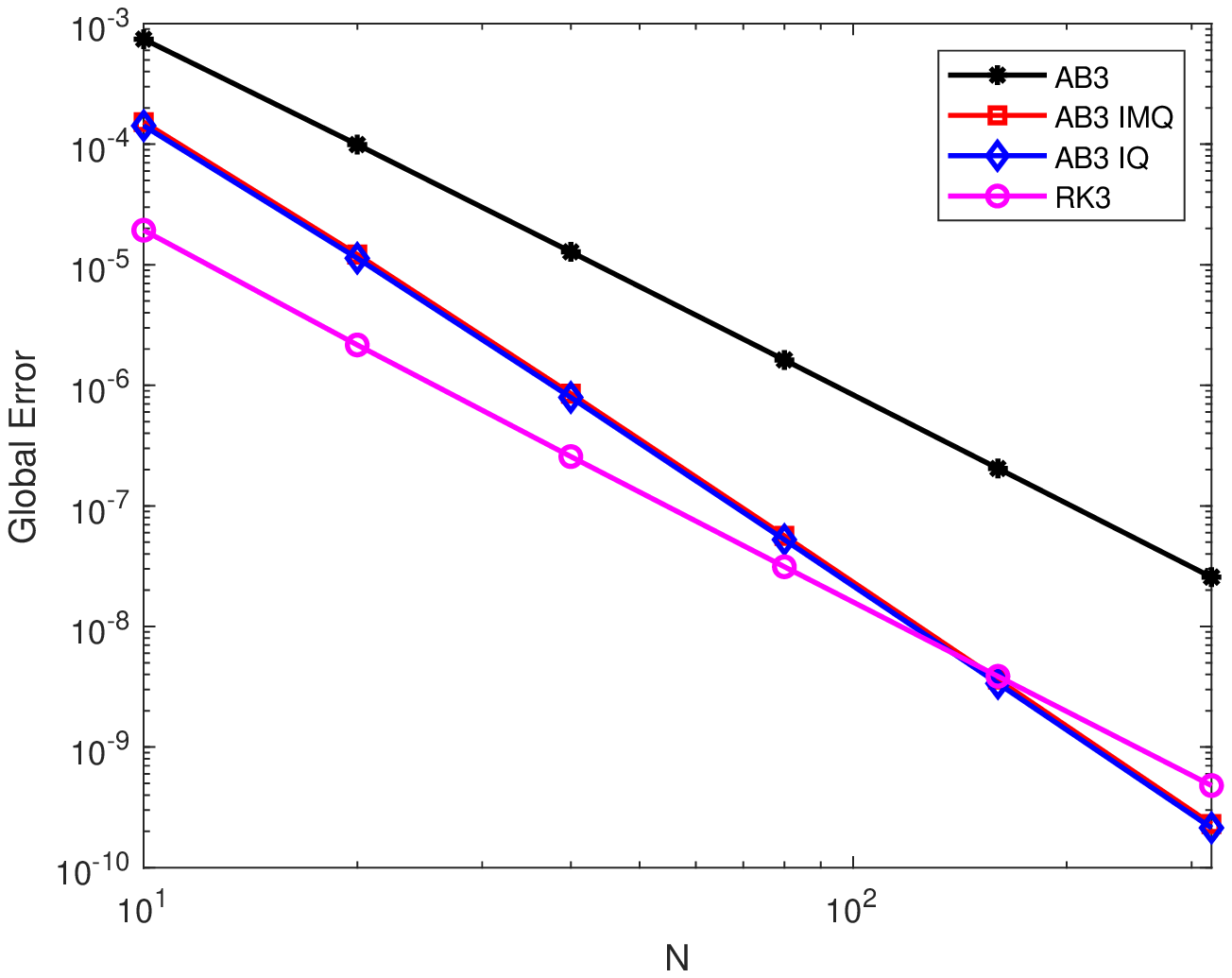} }}
\\
\subfloat{{\includegraphics[trim=1mm  1mm 0.2cm 0cm, clip=true, scale=0.5]{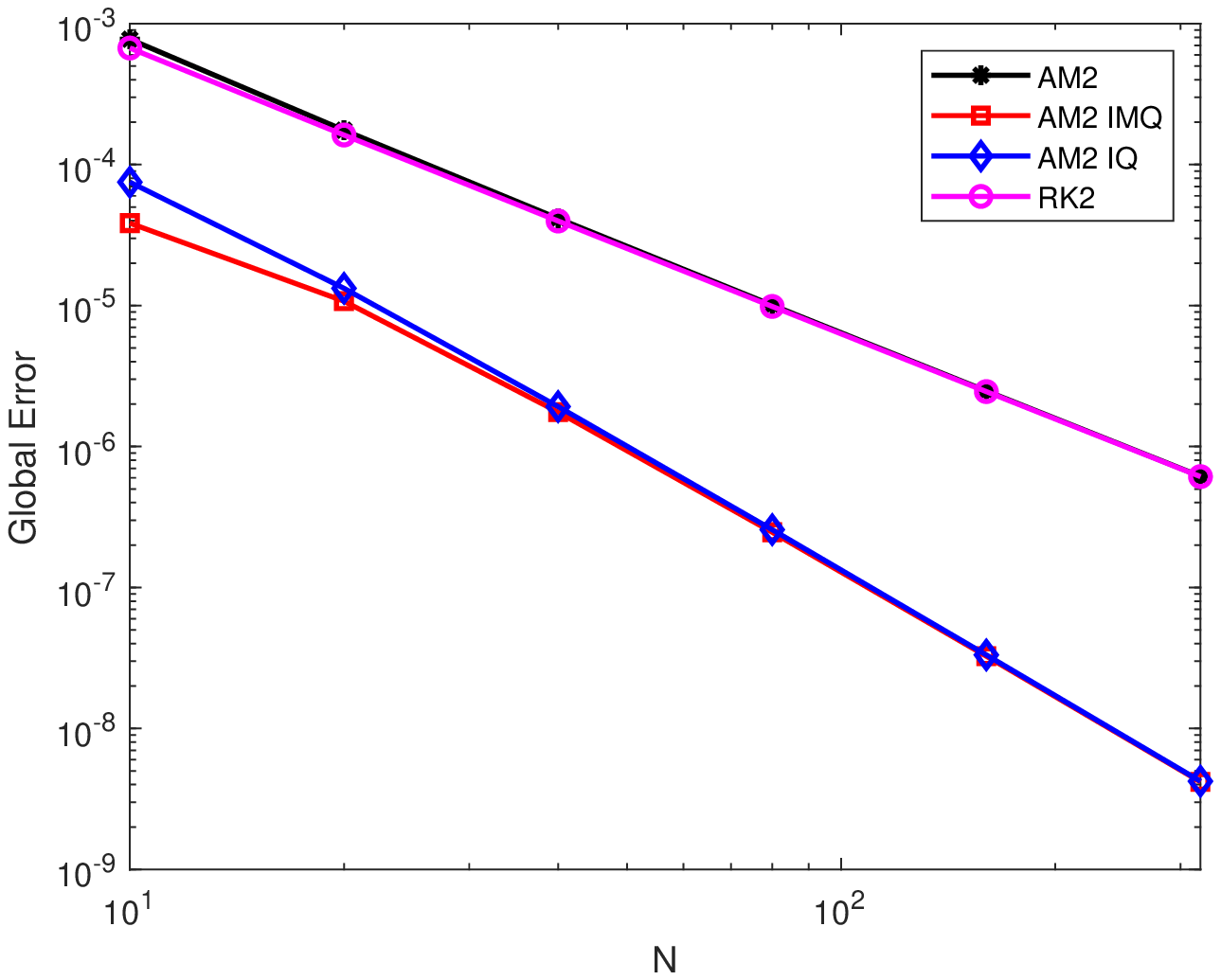} }}
\qquad
\subfloat{{\includegraphics[trim=1mm  1mm 0.2cm 0cm, clip=true, scale=0.5]{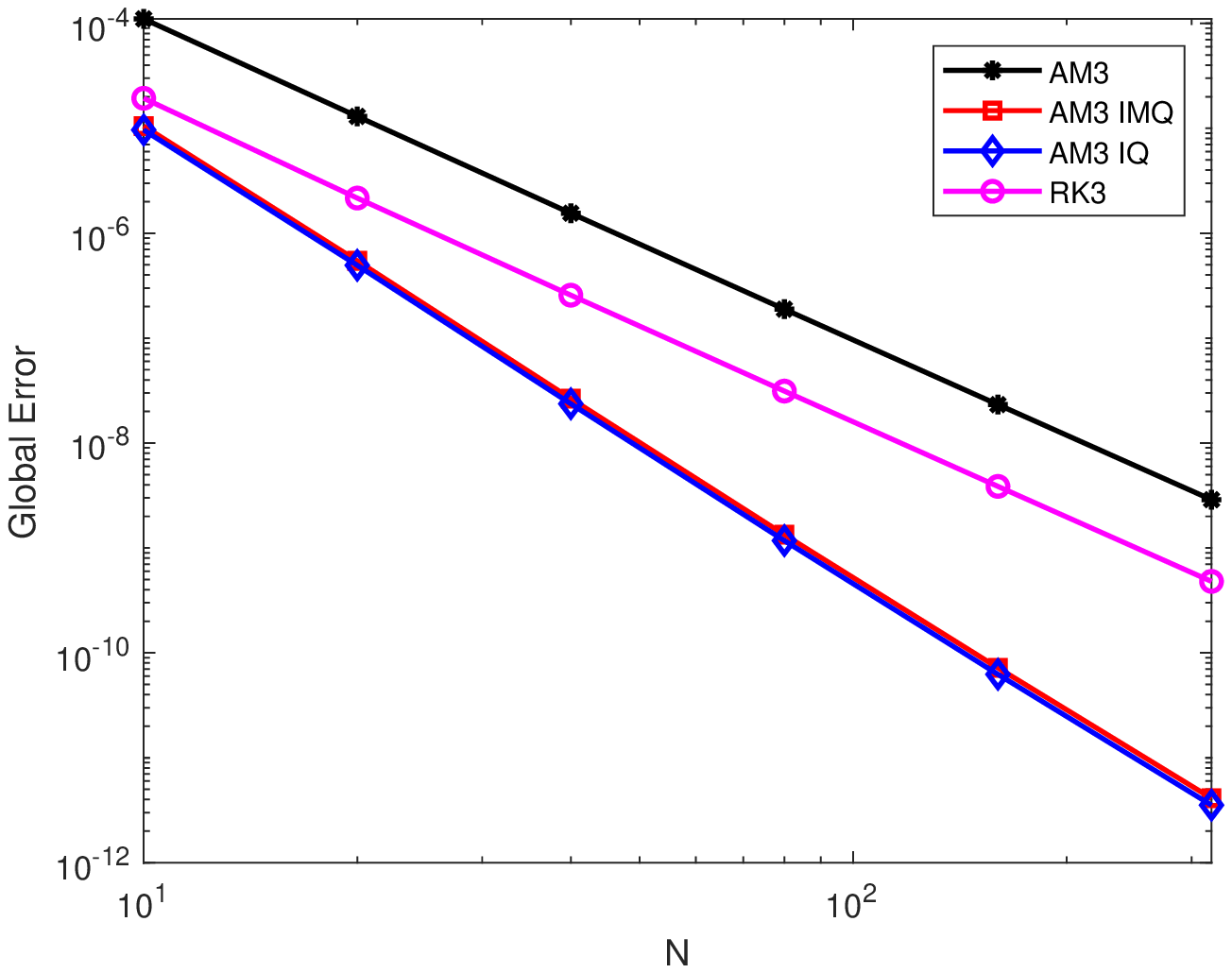} }}

\caption{The global errors versus $N$ in logarithmic scale for example-\ref{Ex1}}%
 \label{fig:7}
\end{figure}

In this section, we present four numerical problems to verify the performance of the proposed methods and compare the results with the original methods to see how these techniques improve the local order of convergence.
\begin{example}\label{Ex1}
{\rm
 We consider the following initial value problem
\begin{equation}
\frac{du}{dt}=-u^2,0<t \leq 1, u(0)=1.
\end{equation}
The exact solution to this problem is $u(t)=\displaystyle\frac{1}{t+1}$. We compute the global errors versus various N and local order of convergence for proposed methods in comparison with original methods. This is shown in Table \ref{Table:1} and \ref{Table:3}  and we plot the same in figure \eqref{fig:7}.  It is observed from Table \ref{Table:1}  upper part that, original RK2 and AB2 methods are second order convergent  whereas the proposed IQ, IMQ-RBF AB2 methods achieves third-order of convergence. The proposed methods also have a better accuracy as compared to conventional RK2 and AB2 methods. The similar behaviour can also be seen in the case of RK3,  AB3,  IQ, IMQ-RBF AB3 methods, see Table \ref{Table:1}  lower part.   In Table \ref{Table:3} upper part, we compare RK2, AM2 with proposed RBF-IMQ and IQ- AM2 methods and concludes that RBF-IMQ and IQ- AM2  achieves third-order convergence with less error whereas RK2, AM2  have second order accuracy only.  Similar behaviour can be seen in the case of RK3, AB3, IMQ-AM3, IQ-AM3 methods in the lower part of the Table \ref{Table:3}. From figure \eqref{fig:7}, it is concluded that the proposed RBF-IMQ and IQ methods perform better in terms of order of convergence and accuracy than the original methods.   In Table \ref{Table:S1}, we show the number of FE, storage space taken, and time-elapsed to perform AB2, AM2, AB3, AM3, RK2, RK3, and proposed methods with the mesh size of $N=100$ for  example \ref{Ex1}. We conclude that our proposed methods  takes almost similar or higher  both in terms of storage and time  taken to get the solution when compares to conventional interpolation methods, while giving a better accuracy with a higher order of convergence.} 
\end{example}

\begin{table}[!h]
\scriptsize
\begin{center}
\caption{ Global errors at $t=1$ and order of convergence for example (\ref{Ex1})}
\begin{tiny}
\label{Table:1}
\begin{tabular}{c| c| c| c| c}
\hline
     \multirow{2}{*}{N} & \multicolumn{4}{c}{Global Error (Order of convergence)}\\\cline{2-5}
    & RK2&AB2  &  RBF IMQ AB2  & RBF IQ AB2 \\
     \hline
10&6.712212827543196e-04(-----) &3.034213293051979e-03(-----)&1.014577580230713e-03(-----)&5.038309526130824e-04(-----)\\
20&1.620903309670352e-04(2.0450)& 7.717017538266813e-04(1.9752)&1.319855656692903e-04(2.9424)&6.437212122212266e-05(2.9684)\\
40&          3.979434794565417e-05(2.0262)& 1.942291044479960e-04(1.9903)&1.670368395367827e-05(2.9821)&8.070554016970100e-06(2.9957)\\
80&9.857160125692488e-06(2.0133)&4.869943879992622e-05(1.9958)&2.097728918681874e-06(2.9933)& 1.008518872303021e-06(3.0004)\\
160&          2.452849796386047e-06(2.0067)&1.219136102492691e-05(1.9980)&2.627413366873554e-07(2.9971)&1.259923502194837e-07(3.0008)\\
320&     6.117820574580435e-07(2.0034)&3.049824816026003e-06(1.9991)&3.287300842647056e-08(2.9987)&1.574287422645426e-08(3.0006)\\[0.5 ex]
\hline
   N & RK3&AB3  &  RBF IMQ AB3  & RBF IQ AB3 \\
     \hline
10&1.933740854498378e-05(-----)                         &7.430910927918588e-04(-----)&1.525244127912884e-04(-----)&1.421782464514632e-04(-----)\\
20&2.162658711668541e-06(3.1605)     & 9.925706387875488e-05(2.9043)& 1.216124766295623e-05(3.6487)&1.132787643659539e-05(3.6498)\\
40               &2.566012992089028e-07(3.0752)     &1.279667578490962e-05(2.9554)&8.526626621430111e-07(3.8342)& 7.939137788848249e-07(3.8348)\\
80        &3.127752390419403e-08(3.0363)       & 1.623800487515759e-06(2.9783)&5.636826039268072e-08(3.9190)& 5.247294543320180e-08(3.9193)\\
160       &3.861620023748458e-09(3.0178)        & 2.044853487093157e-07(2.9893)& 3.622233790689933e-09(3.9599)&3.371532331097882e-09(3.9601)\\
320    &4.797524044697354e-10(3.0088)     & 2.565496837192427e-08(2.9947)&2.295398315865782e-10(3.9801)& 2.136404941843750e-10(3.9801)\\[0.5 ex]
\hline
\end{tabular}
\end{tiny}
\end{center}
\end{table}

\begin{table}[!h]
\scriptsize
\begin{center}
\caption{ Global errors at $t=1$ and order of convergence for example (\ref{Ex1})}
\begin{tiny}
\label{Table:3}
\begin{tabular}{c| c| c| c| c}
\hline
  \multirow{2}{*}{N} & \multicolumn{4}{c}{Global Error (Order of convergence)}\\\cline{2-5}
    &RK2&AM2&   RBF IMQ AM2  & RBF IQ AM2 \\
     \hline
10&6.712212827543196e-04(-----)                         &7.753684341746392e-04(-----)&3.844720863910300e-05(-----)&7.506991550942921e-05(-----)\\
20&1.620903309670352e-04(2.0450)    &1.754079900484484e-04(2.1442)&1.074390241362355e-05(1.8394)&1.325486070935522e-05(2.5017)\\
40&3.979434794565417e-05(2.026)     &4.148050726904273e-05(2.0802)&1.756382688733460e-06(2.6128)&1.920298273949683e-06(2.7871)\\
80&9.857160125692488e-06(2.0133)     &1.006933784897246e-05(2.0425)&2.466792591304667e-07(2.8319)&2.571422232922060e-07(2.9007)\\
160&2.452849796386047e-06(2.0067)     &2.479461800375393e-06(2.0219)&3.257045100291123e-08(2.9210)&3.323120434384208e-08(2.9520)\\
320&6.117820574580435e-07(2.0033) &6.151142203369986e-07(2.0111)&4.181011092896370e-09(2.9616)&4.222521887697894e-09(2.9764)\\[0.5 ex]
\hline
N  &RK3&AM3&   RBF IMQ AM3  & RBF IQ AM3 \\
     \hline
     10&1.933740854498378e-05(-----)                         &1.104840310713895e-04(-----)&1.050535370195060e-05(-----)&9.647365355580639e-06(-----)\\
20&2.162658711668541e-06(3.1605)&1.298901521773477e-05(3.0885)&5.453207364558565e-07(4.2679)&4.945506376641262e-07(4.2859)\\
40          &2.566012992089028e-07(3.0752)& 1.551538363786520e-06(3.0655)& 2.650622332733832e-08(4.3627)& 2.366657114727388e-08(4.3852)\\
80     &3.127752390419403e-08(3.0363)&1.887581442261421e-07(3.0391)&1.335914157074569e-09(4.3104)&1.173951602595480e-09(4.3334)\\
160&3.861620023748458e-09(3.0178)&2.324890380211997e-08(3.0213)&7.196310214396817e-11(4.2144)&6.241751560054354e-11(4.2333)\\
320          &4.797524044697354e-10(3.0088) &2.883801752950887e-09(3.0111)& 4.112044038606655e-12(4.1293)&3.535283177313886e-12(4.1421)\\[0.5 ex]
\hline
\end{tabular}
\end{tiny}
\end{center}
\end{table}

\begin{table}[h!]
\begin{center}
\begin{tabular}{ l | l | l  | p{2.5cm} }
\hline 
\textbf{Method} & \textbf{FE} & \textbf{Memory (Mb) } & \textbf{time(s)} \\
\hline
 \textit{AB2}  &  200 & 0.0923 &0.00140\\
\textit{AM2}&200&0.0932&0.00145  \\
\textit{RK2}&200&0.0940&0.00112\\
\textit{AB3}&300&0.0942&0.00162\\
\textit{AM3}&300&0.0943&0.00167\\
\textit{RK3}&200&0.0937&0.00148\\
 \textit{RBF AB2 IMQ}  & 200&0.0945&0.00160\\
  \textit{RBF AB2 IQ}  & 200&0.0949&0.00140\\
   \textit{RBF AB3 IMQ}  & 300&0.0954&0.00159\\
  \textit{RBF AB3 IQ}  & 300&0.0953&0.00166\\
   \textit{RBF AM2 IMQ}  & 200&0.0947&0.00188\\
  \textit{RBF AM2 IQ}  & 200&0.0954&0.00186\\
   \textit{RBF AM3 IMQ}  & 300&0.0951&0.00192\\
  \textit{RBF AM3 IQ}  & 300&0.0952&0.00191\\
\hline
\end{tabular}
\caption{FE, memory and time-lapsed for solving example \ref{Ex1} with $N=100$}
\label{Table:S1}
\end{center}
\end{table}

\begin{figure}[!htb]
\centering
\subfloat{{\includegraphics[trim=1mm  1mm 0.2cm 0cm, clip=true, scale=0.5]{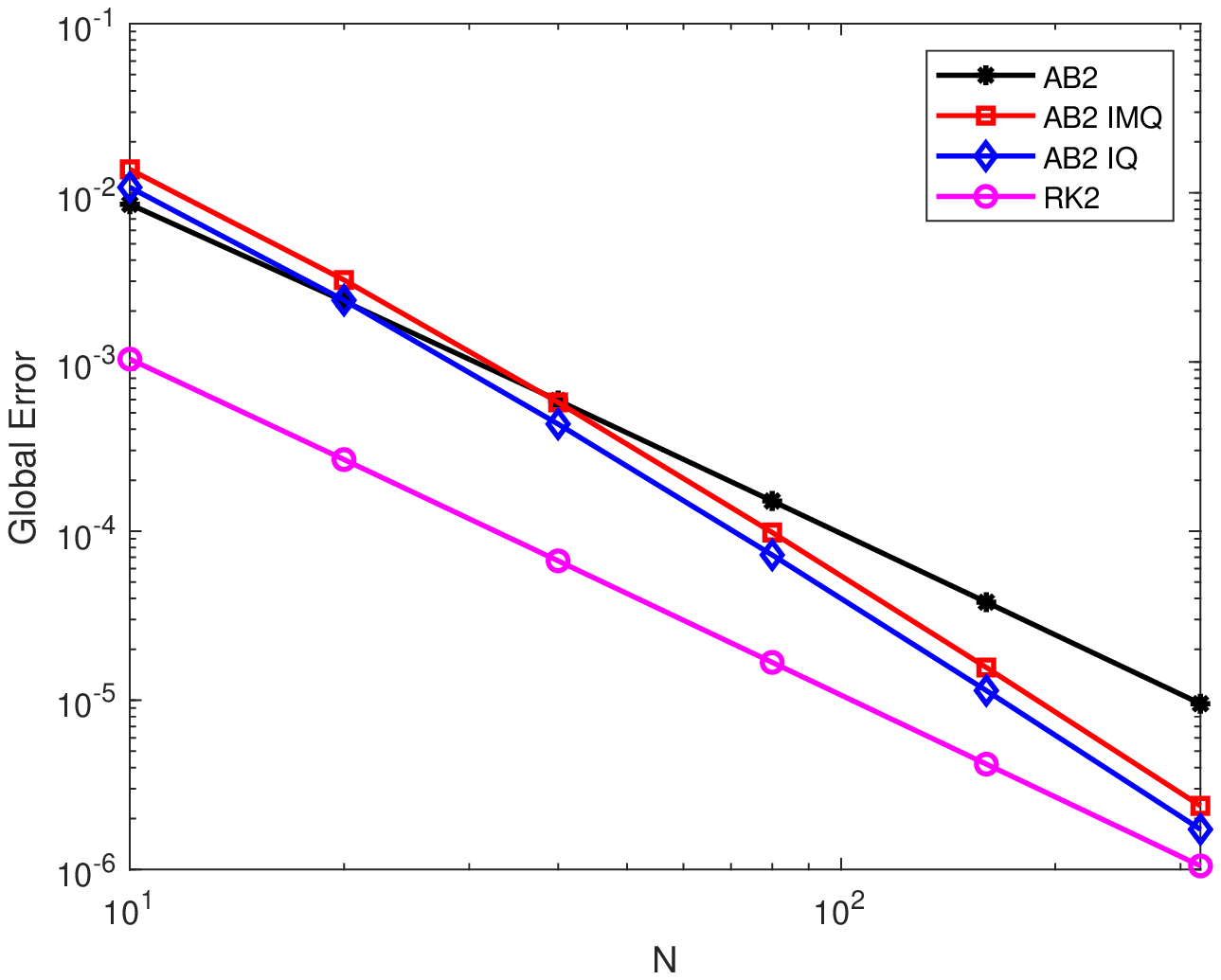} }}
\qquad
\subfloat{{\includegraphics[trim=1mm  1mm 0.2cm 0cm, clip=true, scale=0.5]{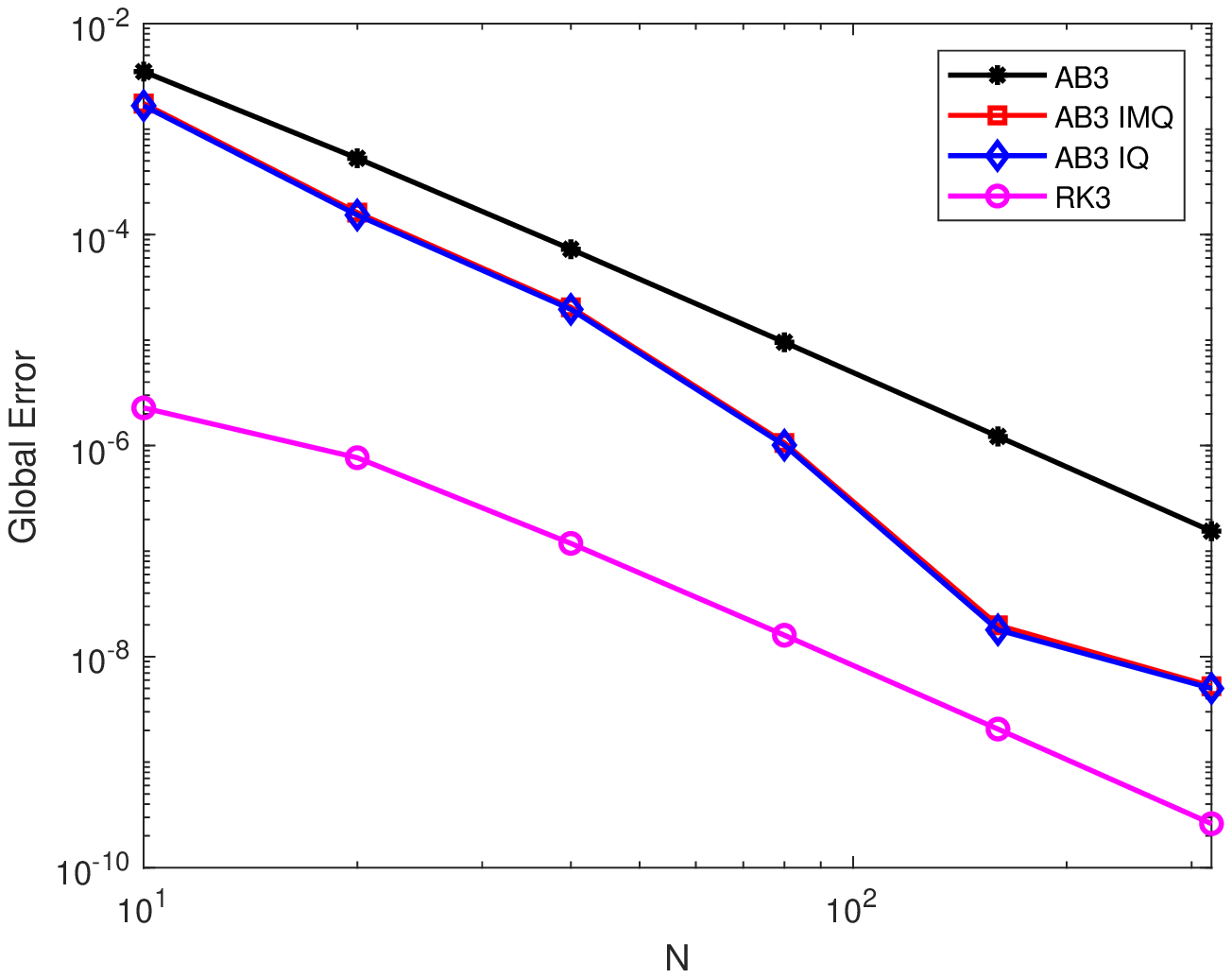} }}
\\
\subfloat{{\includegraphics[trim=1mm  1mm 0.2cm 0cm, clip=true, scale=0.5]{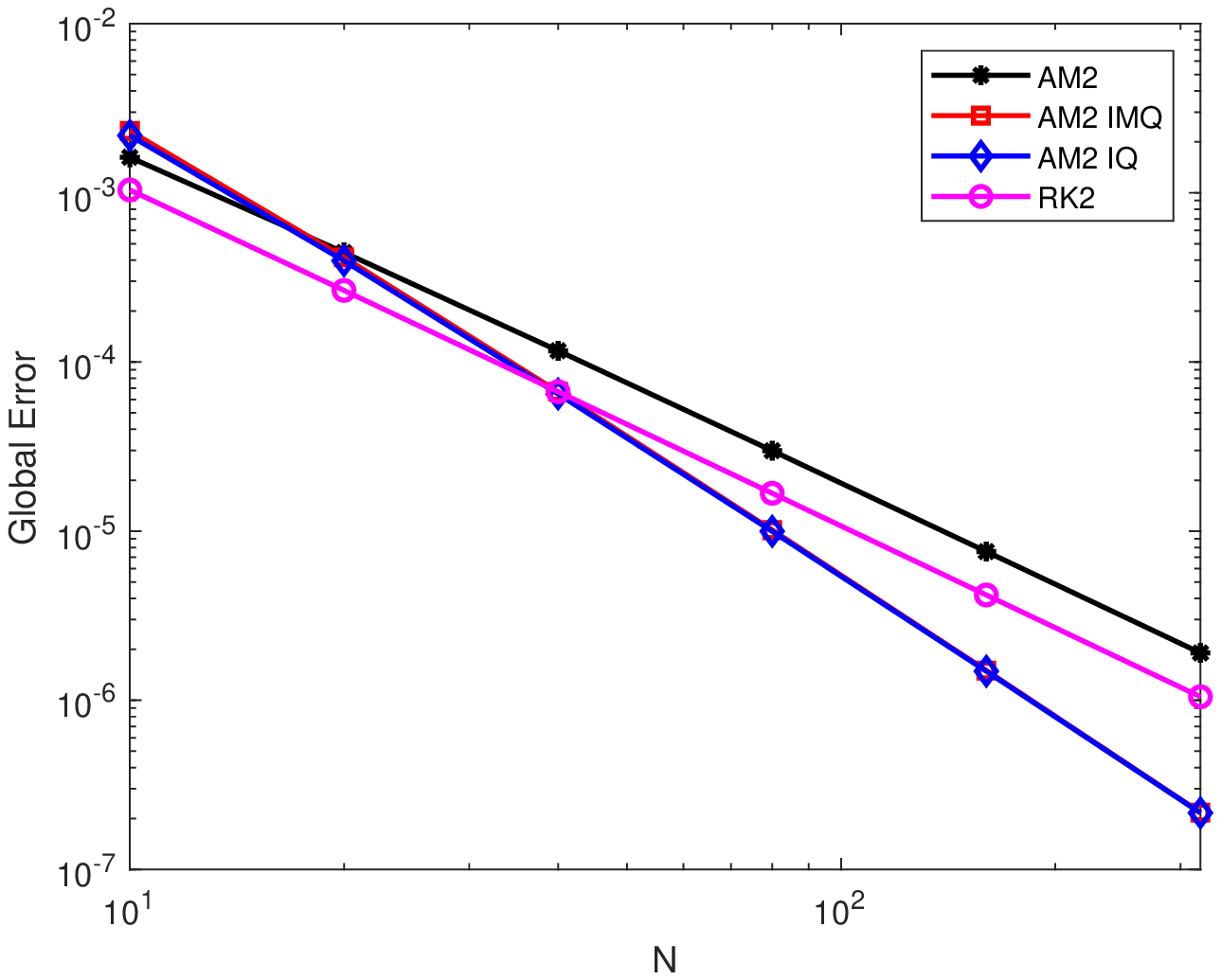} }}
\qquad
\subfloat{{\includegraphics[trim=1mm  1mm 0.2cm 0cm, clip=true, scale=0.5]{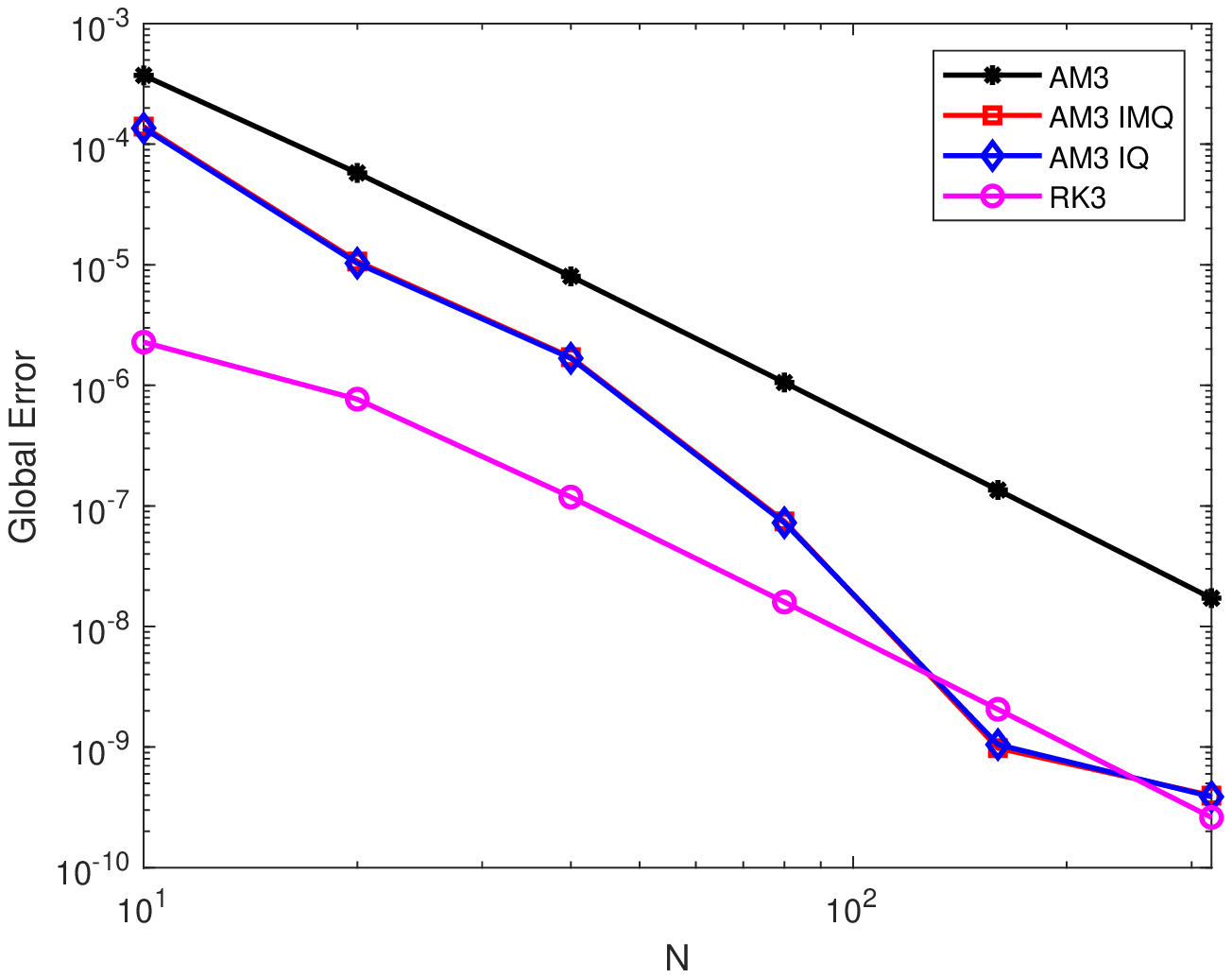} }}
\caption{The global errors versus $N$ in logarithmic scale for example-\ref{Ex2}}%
 \label{fig:9}
\end{figure}

\begin{example}\label{Ex2}
 {\rm
 We consider the following IVP,
\begin{equation}\label{ex2}
u'=\frac{2t^2-u}{t^2u-t},1<t \leq 2, u(1)=2.
\end{equation}
This is a non separable problem, and the exact solution to the differential equation is $u(t)=\frac{1}{t}+\sqrt{\frac{1}{t^2}+4t-4}.$ 
Tables \ref{Table:5} and \ref{Table:7} shows the global errors versus various N and local convergence orders for original and proposed methods respectively.  Figure \ref{fig:9}  show the global errors against N for various methods. In this case also, we observed almost similar behaviour as seen in the earlier case. }
\end{example}

\begin{table}[!h]
\scriptsize
\begin{center}
\caption{ Global errors at $t=2$ and order of convergence for example (\ref{Ex2})}
\begin{tiny}
\label{Table:5}
\begin{tabular}{c| c| c| c| c}
\hline
         \multirow{2}{*}{N} & \multicolumn{4}{c}{Global Error (Order of convergence)}\\\cline{2-5}
    & RK2& AB2  & RBF IMQ AB2 & RBF IQ AB2 \\
     \hline
10&1.039301015929084e-03(-----)                         &8.553947265228956e-03(-----)&1.375772177980528e-02(-----)&1.076861121238171e-02(-----)\\
20&2.648179729285438e-04(1.9725)     &2.286724780864446e-03(    1.9033)&3.053052107214516e-03(     2.1719)&2.316816647421671e-03 (2.2166)\\
40     &6.672648346706112e-05(1.9887)     & 5.919225468034028e-04(     1.9498)& 5.766082861891064e-04(     2.4046)& 4.294793096084604e-04 (2.4315)\\
80&1.673982404692964e-05(1.9950)     &1.506372305621895e-04(     1.9743)&9.818148405393856e-05(     2.5541)& 7.226955744332741e-05     (2.5711)\\
160&4.191773539830024e-06(1.9977)     & 3.800001272402653e-05(     1.9870)&1.562070208782629e-05(     2.6520)& 1.140874991945040e-05 (2.6632)\\
320&1.048766459454953e-06(1.9989)     &9.543155811808646e-06(     1.9935)&2.374327692677980e-06(     2.7179)&1.724804997493834e-06     (2.7256)\\[0.5 ex]
\hline
   N & RK3& AB3  & RBF IMQ AB3 & RBF IQ AB3\\
     \hline
10     &2.278614770201415e-06(-----)                         &3.505315372185347e-03(-----)&1.746638087680719e-03(-----)&1.666341319968456e-03(-----)\\
20&7.675713571408949e-07(1.5698)     &5.294723945774571e-04(2.7269)&1.615198916438132e-04(     3.4348)&1.533183932256321e-04 (3.4421)\\
40     &1.184185203229049e-07(2.6964)     &7.281496346100980e-05(     2.8622)&2.046050758819362e-05(     2.9808)&1.957444451683443e-05     (2.9695)\\
80&1.595202769877346e-08(2.8921)     &9.546965498152815e-06(     2.9311)&1.061039405669817e-06(     4.2693)&1.009575599120183e-06     (4.2772)\\
160&2.057860548632107e-09(2.9545)     &1.222150990631832e-06(     2.9656)&1.982208175377309e-08(     5.7422)&1.793905957114816e-08     (5.8145)\\
320&2.609752414173272e-10(2.9792)     &1.545974748218271e-07(     2.9828)&5.241231892938458e-09(     1.9191)& 4.998635727559986e-09     (1.8435)\\[0.5 ex]
\hline
\end{tabular}
\end{tiny}
\end{center}
\end{table}

\begin{table}[!h]
\scriptsize
\begin{center}
\caption{Global errors at $t=2$ and order of convergence for example (\ref{Ex2})}
\begin{tiny}
\label{Table:7}
\begin{tabular}{c| c| c| c| c}
\hline
     \multirow{2}{*}{N} & \multicolumn{4}{c}{Global Error (Order of convergence)}\\\cline{2-5}
    &RK2 & AM2  & RBF IMQ AM2 & RBF IQ AM2\\
     \hline
10     &1.039301015929084e-03(-----)                         &1.621913472066527e-03(-----)&2.319367796405292e-03(-----)&2.181851550580038e-03(-----)\\
20&2.648179729285438e-04(1.9725)     &     4.438483940938376e-04     (1.8696)&     4.153083431051030e-04(     2.4815)&     3.967344422131092e-04(     2.4593)\\
40&6.672648346706112e-05(1.9887)     &     1.164809592335558e-04     (1.9300)&     6.659715719026238e-05(     2.6407)&     6.470162610971997e-05(     2.6163)\\
80     &1.673982404692964e-05(1.9950)     &     2.987309969038066e-05     (1.9632)&     1.014025616052905e-05(     2.7154)&     9.973896294646067e-06(     2.6976)\\
160     &4.191773539830024e-06(1.9977)     &     7.567079837400570e-06     (1.9810)&     1.498578460878974e-06(     2.7584)&     1.485155865399435e-06(     2.7475)\\
320     &1.048766459454953e-06(1.9989)     &     1.904441726541251e-06     (1.9904)&     2.168098940380503e-07(     2.7891)&     2.157808984826204e-07 (    2.7830)\\[0.5 ex]
\hline
N    &RK3 & AM3  & RBF IMQ AM3 & RBF IQ AM3\\
     \hline
10     &2.278614770201415e-06(-----)                         & 3.716811923353269e-04(-----)&1.397197124024530e-04(-----)&1.357782832358545e-04(-----)\\
20&7.675713571408949e-07(1.5698)     &     5.768342216594391e-05     (2.6878)&1.065042511028125e-05  (   3.7136)& 1.029544897512835e-05((     3.7212)\\
40&1.184185203229049e-07(2.6964) &     8.019406527726147e-06    ( 2.8466)&1.710154137857245e-06  (   2.6387)&1.675849211402181e-06(     2.6190)\\
80     &1.595202769877346e-08(2.8921)  &     1.056329054893723e-06   (  2.9244)&7.442740512786372e-08  (   4.5221)&7.251818168185764e-08(     4.5304)\\
160     &2.057860548632107e-09(2.9545)     &     1.355158145699420e-07 (    2.9625)& 9.796270461492895e-10(     6.2475)&1.048077180598739e-09(     6.1125)\\
320&2.609752414173272e-10(2.9792)   &     1.716002939744499e-08(     2.9813)&3.959388372720696e-10(     1.3070)&3.874256471192439e-10(     1.4358)\\[0.5 ex]
\hline
\end{tabular}
\end{tiny}
\end{center}
\end{table}

\begin{figure}[!htb]
\centering
\subfloat{{\includegraphics[trim=1mm  1mm 0.2cm 0cm, clip=true, scale=0.5]{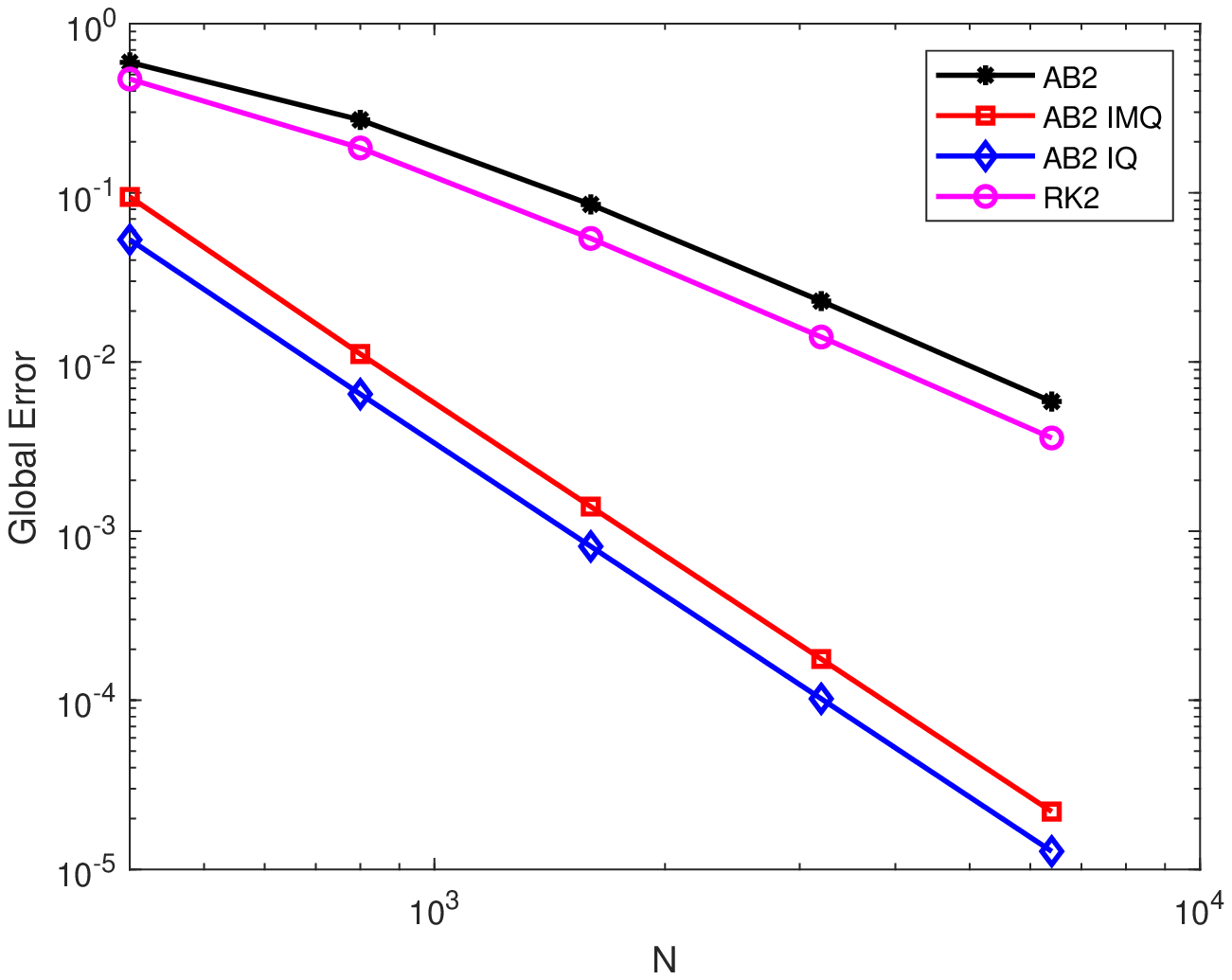} }}
\qquad
\subfloat{{\includegraphics[trim=1mm  1mm 0.2cm 0cm, clip=true, scale=0.5]{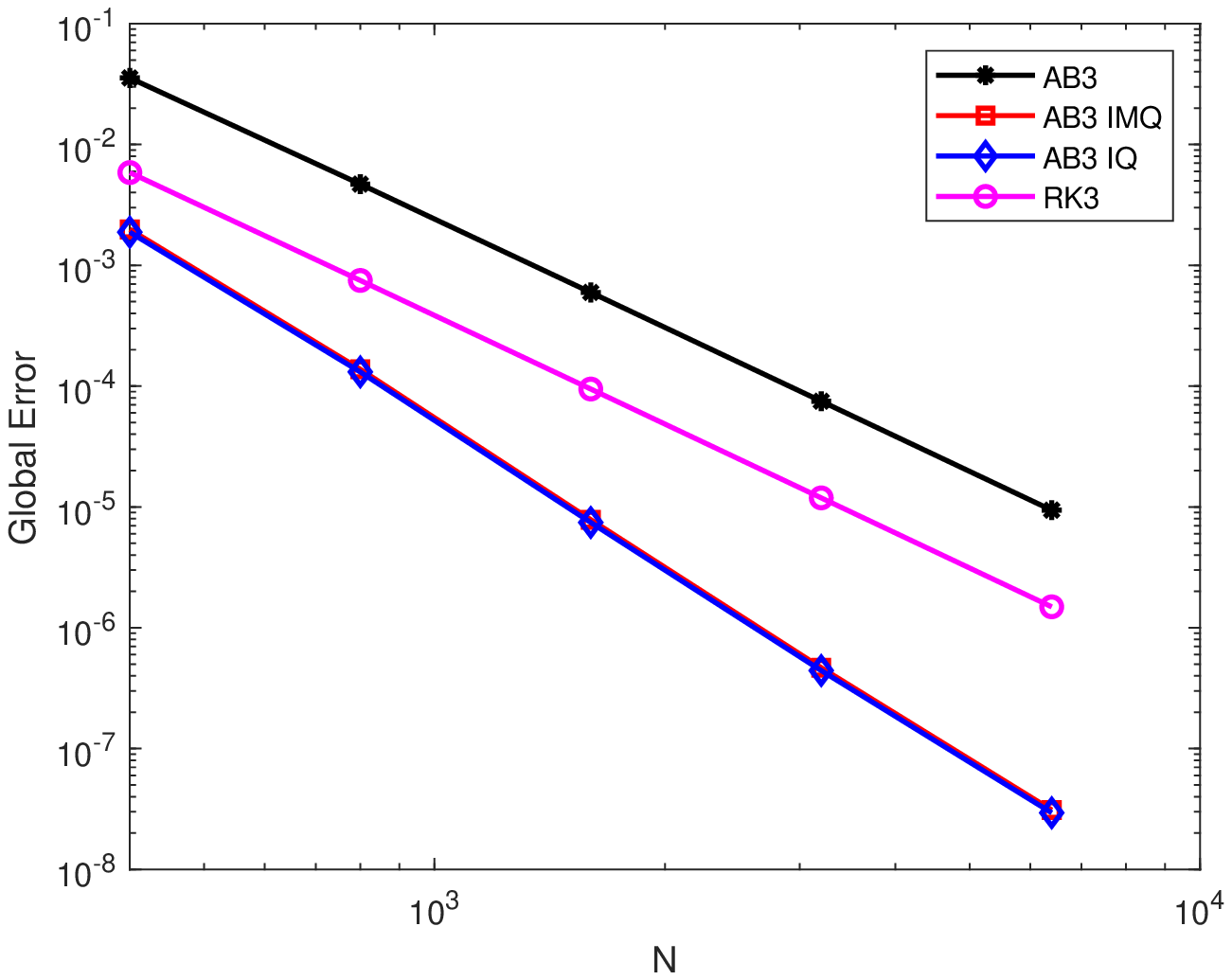} }}
\\
\subfloat{{\includegraphics[trim=1mm  1mm 0.2cm 0cm, clip=true, scale=0.5]{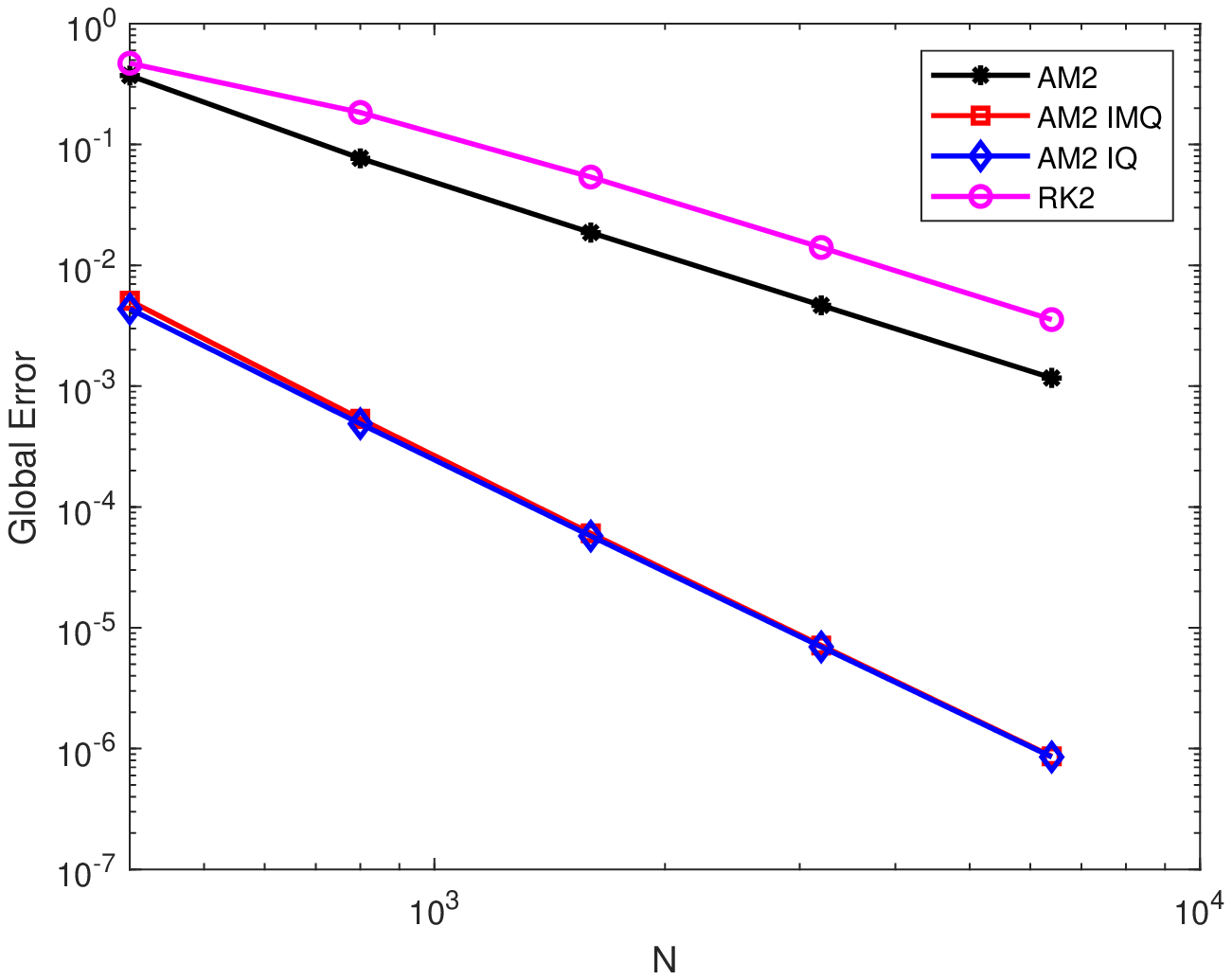} }}
\qquad
\subfloat{{\includegraphics[trim=1mm  1mm 0.2cm 0cm, clip=true, scale=0.5]{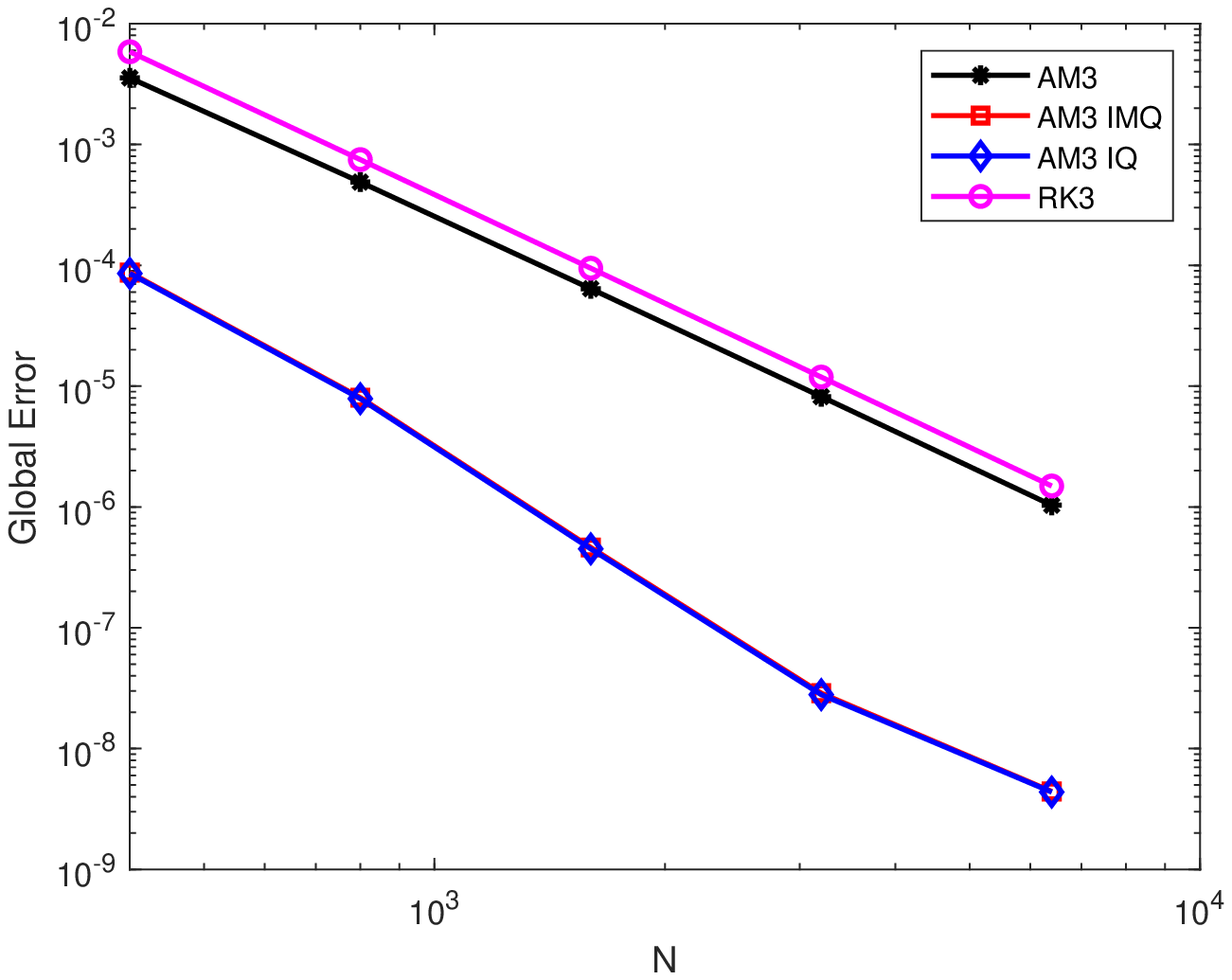} }}
\caption{The global errors versus $N$ in logarithmic scale for example-\ref{Ex3}}
 \label{fig:11}
\end{figure}

\begin{example}\label{Ex3}
 {\rm
Next we consider the following problem
\begin{equation}
u'=-4t^3u^2,-10<t \leq 0,
\end{equation}
with the initial condition $u(-10)=\dfrac{1}{10001}.$ We use this example to verify the results of the proposed methods when applied to a non separable, stiff problem, the solution of which changes rapidly.
The exact solution to the differential equation is $u(t)=\frac{1}{t^4+1}.$  Tables \ref{Table:9} and \ref{Table:11} shows the global errors versus various N and local convergence orders for different proposed methods.  Figure \ref{fig:11} show the global errors against N for different proposed methods. }
\end{example}

\begin{table}[!h]
\scriptsize
\begin{center}
\caption{ Global errors at $t=0$ and order of convergence for example (\ref{Ex3})}
\begin{tiny}
\label{Table:9}
\begin{tabular}{c| c| c| c| c}
\hline
         \multirow{2}{*}{N} & \multicolumn{4}{c}{Global Error (Order of convergence)}\\\cline{2-5}
    &RK2 & AB2  & RBF IMQ AB2 & RBF IQ AB2\\
     \hline
400&4.704254331124198e-01(-----)                         &     5.903133461152865e-01(     0.5170)&9.460567630474581e-02     (4.5277)&5.282492037312436e-02(     3.5906)\\
800&1.841282118830222e-01(1.3533)   &     2.698253771402011e-01(     1.1295)&1.114652209122968e-02     (3.0853)& 6.455035757235272e-03(     3.0327)\\
1600     &5.375587697270057e-02(1.7762)     &     8.536404996574876e-02(     1.6603)&1.395850828363576e-03     (2.9974)&8.119637171259964e-04(     2.9909)\\
3200&1.404512575332451e-02(1.9363)    &     2.289143158128715e-02(     1.8988)&1.753004553235460e-04     (2.9932)&1.020338862329773e-04(     2.9924)\\
6400     &3.553553492288697e-03(1.9827)     &     5.832806784639821e-03(1.9725)&2.197209075882611e-05     (2.9961)& 1.278926515624335e-05(     2.9960)\\[0.5 ex]
\hline
 N   &RK3 & AB3  & RBF IMQ AB3 & RBF IQ AB3\\
     \hline
400&5.854329521159607e-03(-----)                         &     3.551392823076027e-02(     2.6144)&1.979811601221737e-03(     4.0225)&1.879952358237125e-03(     4.0480)\\
800&7.493169352092988e-04(2.9659)    &    4.694922501482846e-03(     2.9192)&1.376981990173221e-04(     3.8458)&1.310069688484905e-04(     3.8430)\\
1600&9.460311993314541e-05(2.9856)     &     5.961251400129486e-04(     2.9774)&7.834765070935390e-06(     4.1355)&7.438336088405606e-06(     4.1385)\\
3200&1.188305366095488e-05(2.9930)     &     7.498011945050731e-05(     2.9910)&4.666861030955261e-07(     4.0694)&4.425260993201618e-07(     4.0711)\\
6400&1.490943683957191e-06(2.9946)     &     9.401686092758155e-06(     2.9955)&3.102741330529568e-08(     3.9108)&2.941590926397453e-08(     3.9111)\\[0.5 ex]
\hline
\end{tabular}
\end{tiny}
\end{center}
\end{table}

\begin{table}[!h]
\scriptsize
\begin{center}
\caption{ Global errors at $t=0$ and order of convergence for example (\ref{Ex3})}
\begin{tiny}
\label{Table:11}
\begin{tabular}{c| c| c| c| c}
\hline
    \multirow{2}{*}{N} & \multicolumn{4}{c}{Global Error (Order of convergence)}\\\cline{2-5}
    &RK2 & AM2  & RBF IMQ AM2 & RBF IQ AM2\\
     \hline
400&4.704254331124198e-01(-----)                         &     3.720550144943622e-01(     6.9022)&5.087211829438054e-03(     3.4221)&4.350402619599736e-03(     3.2880)\\
800&1.841282118830222e-01(1.3533)     &     7.685743832650860e-02(     2.2753)&5.355129181516816e-04(     3.2479)&4.881593670340578e-04(     3.1557)\\
1600&5.375587697270057e-02(1.7762)&1.865512210424147e-02(2.0426)&6.045050353797876e-05(     3.1471)&5.743808932368033e-05(3.0873)\\
3200     &1.404512575332451e-02(1.9363)     &     4.660509234408083e-03(     2.0010)&7.139430071578801e-06(     3.0819)&6.949383444254664e-06(     3.0472)\\
6400&3.553553492288697e-03(1.9827)     &     1.168767830167283e-03(     1.9955)&8.640540185567147e-07(     3.0466)&8.521166778852063e-07(     3.0278)\\[0.5 ex]
\hline
  N  &RK3 & AM3  & RBF IMQ AM3 & RBF IQ AM3\\
     \hline

400&5.854329521159607e-03(-----) &     3.554794861643673e-03(     2.7272)&8.701242429842715e-05(     2.8951)&8.556512650192971e-05(     2.9329)\\
800&7.493169352092988e-04(2.9659)     &     4.884174912045403e-04(     2.8636)&8.029627526107141e-06(     3.4378)&7.868121421283902e-06(     3.4429)\\
1600&9.460311993314541e-05(2.9856)     &     6.398864657386483e-05(     2.9322)&4.615079721714821e-07(     4.1209)&4.498488695992009e-07(     4.1285)\\
3200&1.188305366095488e-05(2.9930)     &     8.186059585302274e-06(     2.9666)& 2.872902371553465e-08(     4.0058)& 2.795302189717574e-08(     4.0084)\\
6400&1.490943683957191e-06(2.9946)     &     1.033128076421974e-06(     2.9862)&4.415429355475453e-09(     2.7019)&4.367913808422941e-09(     2.6780)\\[0.5 ex]
\hline
\end{tabular}
\end{tiny}
\end{center}
\end{table}

\begin{figure}[!htb]
\centering
\subfloat{{\includegraphics[trim=1mm  1mm 0.2cm 0cm, clip=true, scale=0.5]{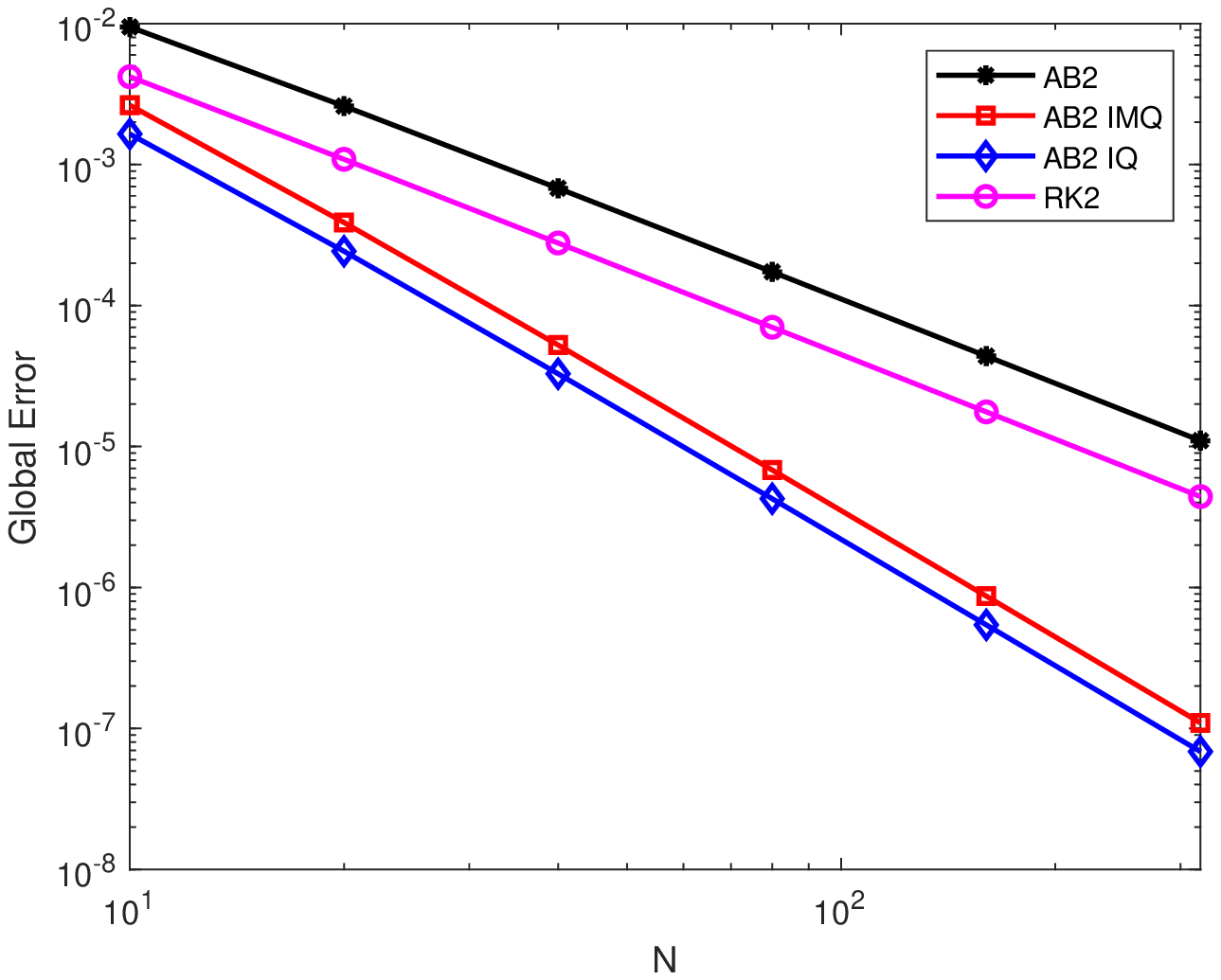} }}
\qquad
\subfloat{{\includegraphics[trim=1mm  1mm 0.2cm 0cm, clip=true, scale=0.5]{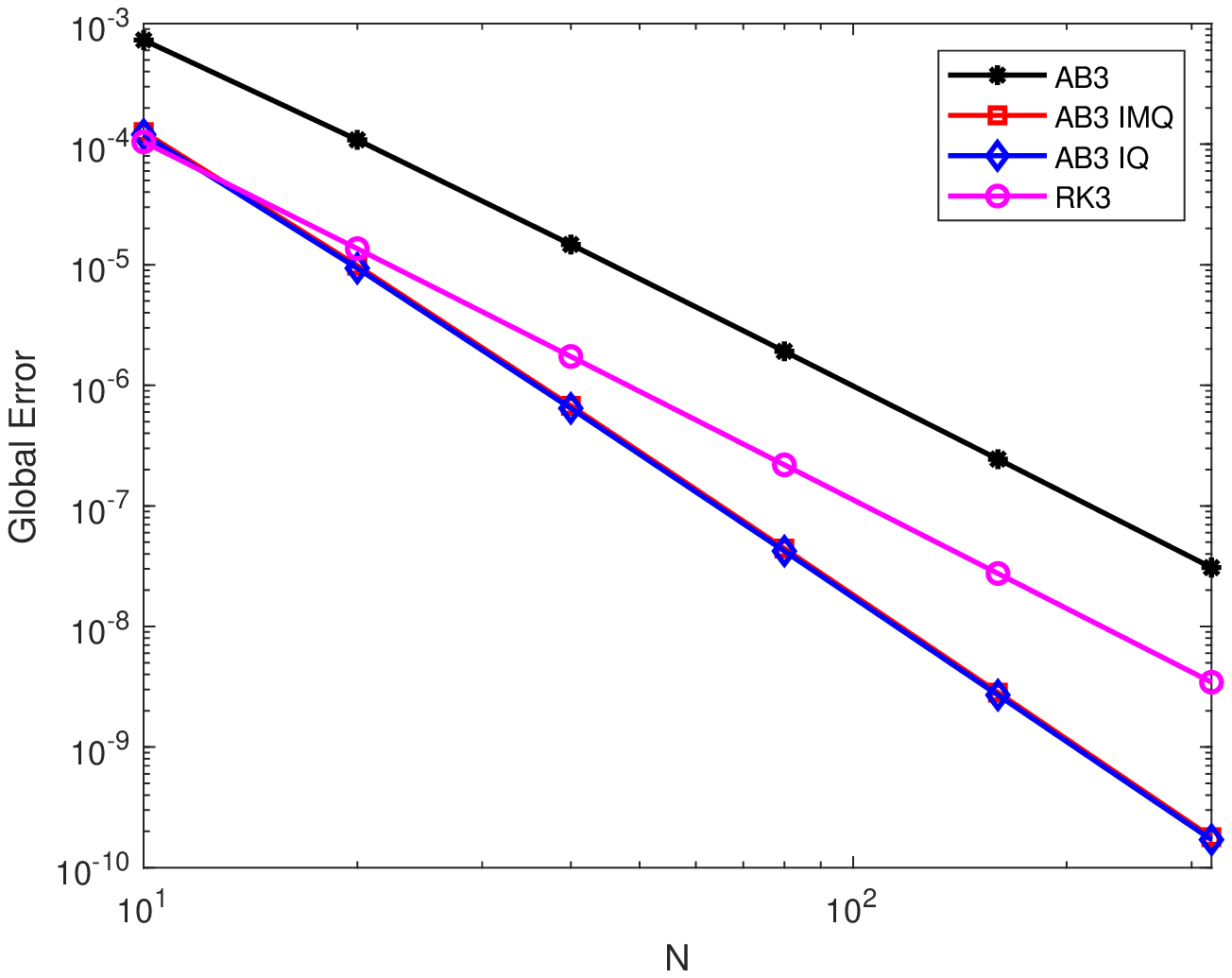} }}
\\
\subfloat{{\includegraphics[trim=1mm  1mm 0.2cm 0cm, clip=true, scale=0.5]{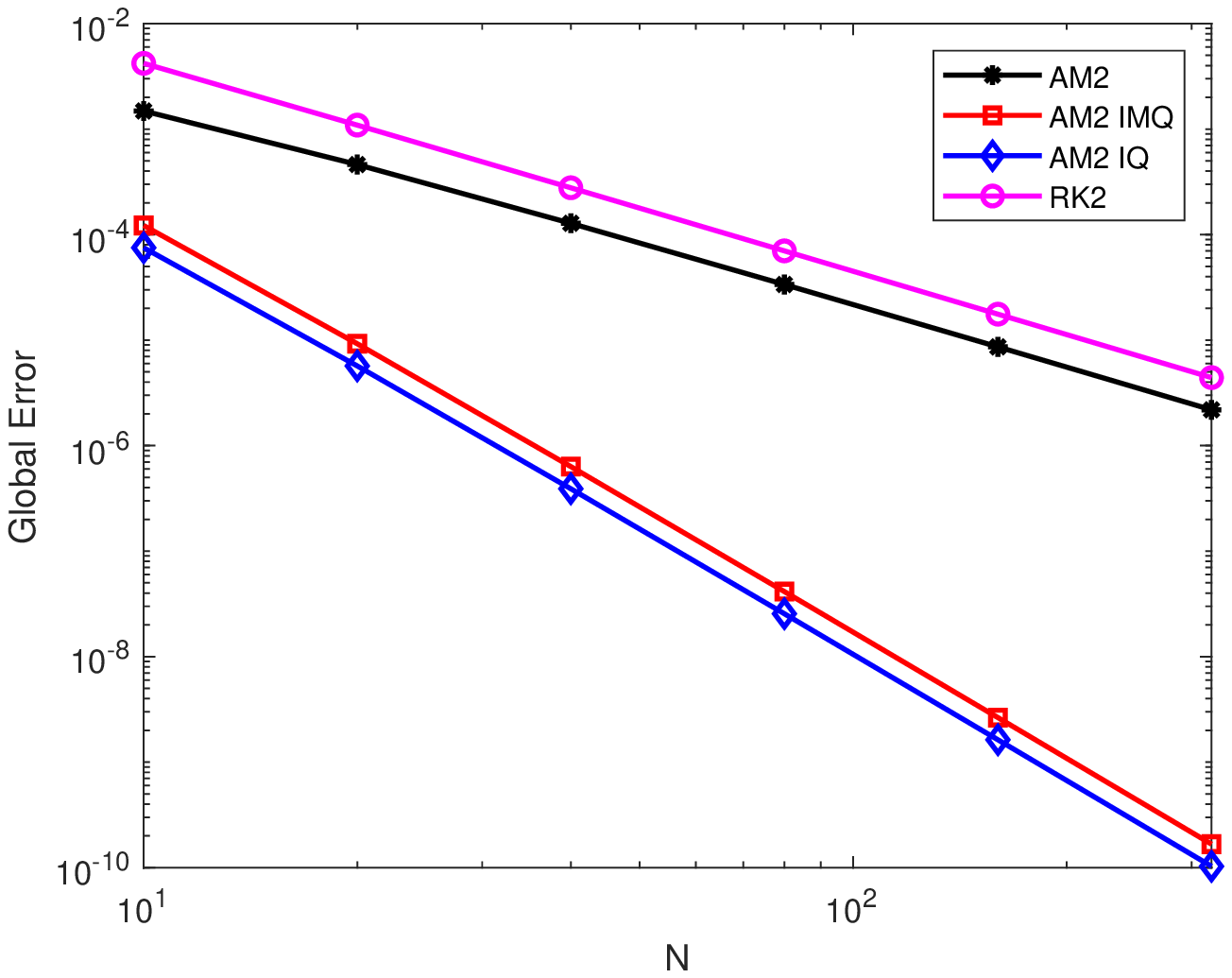} }}
\qquad
\subfloat{{\includegraphics[trim=1mm  1mm 0.2cm 0cm, clip=true, scale=0.5]{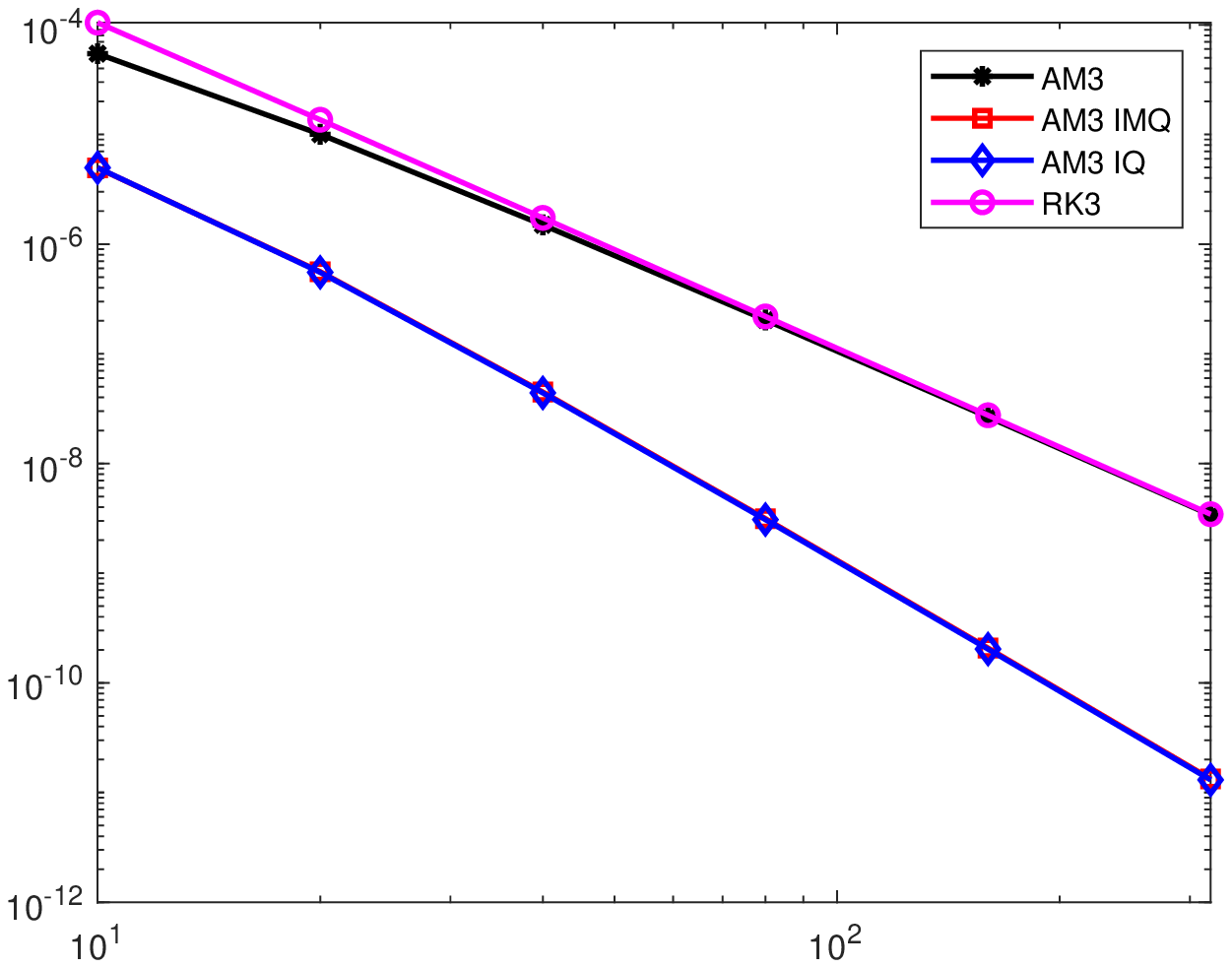} }}
\caption{The global errors versus $N$ in logarithmic scale for example-\ref{Ex4}}
 \label{fig:13}
\end{figure}

\begin{example}\label{Ex4}
 {\rm
Finally, we consider the following initial value problem
\begin{equation}
\frac{du}{dt}=u+2,0<t \leq 1,
\end{equation}
with initial condition
\begin{equation}
u(0)=-1.
\end{equation}
The exact solution to this differential equation is $u(t)=e^t-2$. Tables \ref{Table:13} and \ref{Table:15} shows the global errors versus various N and local convergence orders for different proposed methods.   Figure \ref{fig:13}  shows the global errors against N for different proposed methods. 
}
\end{example}

\begin{table}[!h]
\scriptsize
\begin{center}
\caption{ Global errors at $t=1$ and order of convergence for example (\ref{Ex4})}
\begin{tiny}
\label{Table:13}
\begin{tabular}{c| c| c| c| c}
\hline
        \multirow{2}{*}{N} & \multicolumn{4}{c}{Global Error (Order of convergence)}\\\cline{2-5}
    &RK2 & AB2  & RBF IMQ AB2 & RBF IQ AB2\\
     \hline
10&4.200981850820629e-03(-----)  &9.467968119577841e-03(-----)&2.638993106393817e-03(-----)&1.648319236535523e-03(-----)\\
20&1.090774104160142e-03(1.9454)     &     2.601012298556360e-03(     1.8640)&3.884608567352732e-04(     2.7641)&2.428715386061375e-04(     2.7627)\\
40&2.778840880685030e-04(1.9728)     &     6.792824867547509e-04(     1.9370)&5.250103044651677e-05(     2.8874)&3.282183176667175e-05(     2.8875)\\
80&7.012735968714434e-05(1.9864)     &     1.734128452398620e-04(     1.9698)&6.817907165834747e-06(     2.9449)&4.261807716732235e-06(     2.9451)\\
160     &1.761434225810987e-05(1.9932)     &     4.379920178887442e-05(     1.9852)&8.684610666165327e-07(     2.9728)&5.428282044972477e-07(     2.9729)\\
320&4.413926785629130e-06(1.9966)     &     1.100533493658684e-05(     1.9927)&1.095797216299488e-07(     2.9865)&6.848986400243007e-08(     2.9865)\\[0.5 ex]
\hline
   N  &RK3 & AB3  & RBF IMQ AB3 & RBF IQ AB3\\
     \hline
10&1.045659774348007e-04(-----)                         &7.307854535646419e-04(-----)&1.261298665118371e-04(-----)&1.202909238573469e-04(-----)\\
20     &1.360300818820104e-05(2.9424)     &     1.086626000101898e-04(     2.7496)&9.815480212060912e-06(     3.6837)&9.361590633583106e-06(     3.6836)\\
40     &1.734685968757255e-06(2.9712)     &     1.473078638514203e-05(     2.8829)&6.761295517909716e-07(     3.8597)&6.448734821695723e-07(     3.8597)\\
80     &2.190136697555189e-07(2.9856)     &     1.915342045410995e-06(     2.9432)&4.425224242599057e-08(     3.9335)&4.220683647027101e-08(     3.9335)\\
160     &2.751389294175510e-08(2.9928)     &     2.441164188571321e-07(     2.9720)&2.828651202158028e-09(    3.9676)&2.697916556826385e-09(     3.9676)\\
320&3.447842567005921e-09(2.9964)     &     3.081060218068643e-08(     2.9861)&1.787646697337664e-10(    3.9840)&1.705020569175986e-10(     3.9840)\\[0.5 ex]
\hline
\end{tabular}
\end{tiny}
\end{center}
\end{table}

\begin{table}[!h]
\scriptsize
\begin{center}
\caption{ Global errors at $t=1$ and order of convergence for example (\ref{Ex4})}
\begin{tiny}
\label{Table:15}
\begin{tabular}{c| c| c| c| c}
\hline
    \multirow{2}{*}{N} & \multicolumn{4}{c}{Global Error (Order of convergence)}\\\cline{2-5}
    &RK2 & AM2  & RBF IMQ AM2 & RBF IQ AM2\\
     \hline
10      &4.200981850820629e-03(-----)                                                     & 1.485946965479878e-03(-----)&1.219551075468672e-04(-----)&7.512354998162607e-05(-----)\\
20      &1.090774104160142e-03(1.9454)         &     4.613665446593362e-04(     1.6874)&9.239788074211219e-06(     3.7223)&5.692920712729865e-06(     3.7220)\\
40     &2.778840880685030e-04(1.9728)       &     1.279529299152982e-04(     1.8503)&6.328599858607120e-07(     3.8679)&3.899504021909195e-07(     3.8678)\\
80      &7.012735968714434e-05(1.9864)     &     3.365836437829728e-05(     1.9266)&4.136389708531141e-08(     3.9354)&2.548797151114002e-08(     3.9354)\\
160      &1.761434225810987e-05(1.9932)            &     8.629485029554296e-06(     1.9636)&2.643091190357438e-09(     3.9681)&1.628665202701995e-09(     3.9681)\\
320      &4.413926785629130e-06(1.9966)       &     2.184624928958101e-06(     1.9819)&1.670212856907938e-10(     3.9841)&1.029193397172889e-10(     3.9841)\\[0.5 ex]
\hline
    \multirow{2}{*}{N} & \multicolumn{4}{c}{Global Error (Order of convergence)}\\\cline{2-5}
    &RK3 & AM3  & RBF IMQ AM3 & RBF IQ AM3\\
     \hline
10&1.045659774348007e-04(-----) &2.353885010877499e-02(-----)&6.472686308417464e-04(-----)&6.534085433826942e-04(-----)\\
20&1.360300818820104e-05(2.9424)&     3.554794861643673e-03(     2.7272)&8.701242429842715e-05(     2.8951)&8.556512650192971e-05(     2.9329)\\
40&1.734685968757255e-06(2.9712)&        4.884174912045403e-04(     2.8636)&8.029627526107141e-06(     3.4378)&7.868121421283902e-06(     3.4429)\\
80&2.190136697555189e-07(2.9856) &       6.398864657386483e-05(     2.9322)&4.615079721714821e-07(     4.1209)&4.498488695992009e-07(     4.1285)\\
160  &2.751389294175510e-08(2.9928) &    8.186059585302274e-06(     2.9666)&2.872902371553465e-08(     4.0058)&2.795302189717574e-08(     4.0084)\\
320&3.447842567005921e-09(2.9964)  &     1.033128076421974e-06(     2.9862)&4.415429355475453e-09(     2.7019)&4.367913808422941e-09(     2.6780)\\[0.5 ex]
\hline
\end{tabular}
\end{tiny}
\end{center}
\end{table}

\section{Conclusion}
In this work, we present four improved Adam-Bashforth and four Adam-Moulton methods to tackle initial value problems. Interpolation of inverse-quadratic and inverse multi-quadratic radial basis functions is used in these approaches to find an approximation. These enhanced Adam-Bashforth and Adam-Moulton methods give more accurate results than the original Adam-Moutlon and Adam-Bashforth ODE solver by optimising the free shape parameter of the RBF functions. The consistency, stability regions  and local order of convergence has been investigated. This is a preliminary study on the shape parameter, that is shape parameter is locally identical for the IQ, IMQ-RBF functions. Further study is going on in two directions (i) if the shape parameter is chosen arbitrary, (ii) extension of these methods on non-uniform mesh and these will be reported separately. 

\section*{Acknowledgements}
Samala Rathan has been supported by IIPE, Visakhapatnam, India, under the IRG grant number IIPE/DORD/IRG/001.



\end{document}